\documentclass{amsart}
\usepackage{amssymb}
\usepackage{amsfonts}
\usepackage{amsmath}
\usepackage{lastpage}
\usepackage[legalpaper,bookmarks=true,colorlinks=true,linkcolor=blue,citecolor=blue]{hyperref}
\usepackage{graphicx}%
\usepackage{fancyhdr}
\usepackage{color}
\usepackage[mathlines]{lineno}
\usepackage{lscape}
\usepackage{epsfig}
\usepackage{geometry}
\usepackage{tgbonum}
\fontfamily{qcr}\selectfont

\usepackage[]{listings}

\fontfamily{qcr}\selectfont

\theoremstyle{plain}

\newtheorem{corollary}{Corollary}

\newtheorem{proposition}{Proposition}

\numberwithin{equation}{section}

\newcommand{\Bin}{\bigskip \noindent}

\newcommand{\Ni}{\noindent}

\begin{document}
\Large
\title[Joint limits laws of lower and upper records]{Simultaneous Joint Lower and Upper record values Probability Laws for Absolutely Continuous or Discrete Data}
\author{Gane Samb Lo $^{\dag}$}
\author{Mohammad Ahsanullah $^{\dag\dag}$}
\author{Aladji Babacar Niang $^{\dag\dag\dag}$}

\begin{abstract} This paper investigates the probability density function (\textit{pdf}) of the $(2n-1)$-vector  ($n\geq 1$) of both lower and upper record values for  a sequence of independent random variables with common \textit{pdf} $f$ defined on the same probability space, provided that the lower and upper record times are finite up to $n$. A lot is known about the lower or the upper record values when they are studied separately. When put together, the challenges are a far bigger complicated. The rare results in the literature still present important flaws. This paper begins a new and complete investigation with a few number of records: 2 and 3. Lessons from these simple cases will allow addressing the general formulation of simultaneous joint lower-upper records.\\ 

\noindent $^{\dag}$ Gane Samb Lo.\\
LERSTAD, Gaston Berger University, Saint-Louis, S\'en\'egal (main affiliation).\newline
LSTA, Pierre and Marie Curie University, Paris VI, France.\newline
AUST - African University of Sciences and Technology, Abuja, Nigeria\\
Imhotep Mathematical Center (IMC)\\
gane-samb.lo@edu.ugb.sn, gslo@aust.edu.ng, ganesamblo@ganesamblo.net\\
Permanent address : 1178 Evanston Dr NW T3P 0J9,Calgary, Alberta, Canada.\\

\noindent $^{\dag\dag}$ Mohammad Ahsanullah\\
Department of Management Sciences. Rider University\\
Lawrenceville, New Jersey, USA\\
Email : ahsan@rider.edu\\

\noindent $^{\dag\dag\dag}$ Aladji Babacar Niang\\
LERSTAD, Gaston Berger University, Saint-Louis, S\'en\'egal.\\
Imhotep Mathematical Center (IMC)\\
Email: niang.aladji-babacar@ugb.edu.sn, aladjibacar93@gmail.com\\

\end{abstract}
\maketitle
\newpage
\tableofcontents

\section{Introduction} \label{sec_01_ulrecords}

\noindent Let us $Y_1$, $Y_2$, $\cdots$ be a sequence of real random variables defined on the same probability space $\left(\Omega, \mathcal{A}, \mathbb{P}\right)$. The theory of records deals with the probability laws of the strong or weak upper record values and the associated record times and the lower versions of such mathematical objects. Next, asymptotic results are drawn and applications made in real-life situations.\\
  
\Ni A great deal of that theory is known for \textit{iid} data and for a few cases with dependent data, but usually in the stationary frame (see \cite{nevzorov}). The books by \cite{nevzorov}, \cite{ahsan88}, \cite{ahsan95}, \cite{ahsan04}, \cite{arnold}, etc. introduce to that theory. Recently, in \cite{ahsanloGL2019}, a general formulation of probability laws regardless of the dependence has been given.\\

\Ni Also, a stochastic process view, via the extremal process, is provided as an extension of record value theory (\cite{dwass}, \cite{deheuvels}, \cite{resnickPOS}, etc.)\\

\Ni It happens that results for upper records are easily transferred to lower records by the opposite transform or by the inverse transform for positive data. As a result, most of the results have been given for upper records.\\

\Ni Up to our knowledge, a simultaneous study of upper and lower records is not well documented yet and some of the available results present flaws. However, probability laws of joint lower-upper records are important in a number of situations, for example for studying the discrepancy between the upper records and lower records. For example, the empirical range of a real-valued distribution is an example of difference of record values.\\

\Ni In the beginning of that effort to establish the general law of strong simultaneous lower-upper record values and times, both for absolutely continuous and 
discrete records, we see that this task is far from easy. So, to layout the way, we completely determine the joint simultaneous lower-upper record values laws up to three lower records and upper records as a record process.\\

\Ni Let us be more precise with the adequate notation.\\

\Ni \textbf{Strong upper record times}. Let us put $u(1)=1$ as the first strong upper record time. For any $n\geq 2$, we define, by induction, whenever the $(n-1)$-\textit{th} upper record time $u(n-1)$ exists,

$$
U_{n}=\left\{j>u(n-1), \ Y_j>Y_{u(n-1)}\right\}.
$$

\Bin Hence, for $n\geq 2$, the $n$-\textit{th} upper record time is $u(n)=+\infty$ if $U_n$ is empty and, otherwise

$$
u(n)=\inf U_n.
$$

\Bin  \textbf{Strong lower record times}. Let us put $\ell(1)=1$ as the first strong lower record time. For any $n\geq 2$, we define, by induction, whenever the $(n-1)$-\textit{th} lower record time $\ell(n-1)$ exists,

$$
L_{n}=\left\{j>\ell(n-1), \ Y_j<Y_{\ell(n-1)}\right\}.
$$

\Bin  Hence, for $n\geq 2$, the $n$-\textit{th} lower record time is $\ell(n)=+\infty$ if $L_n$ is empty and, otherwise

$$
\ell(n)=\inf L_n.
$$

\Bin  \textbf{Strong record values}. For each $n\geq 1$ such that $u(n)$ is finite, we have a sequence of  strong upper record values  
$$(Y^{(k)}=Y_{u(k)}, \ 1\leq k\leq n).$$

\Bin For each $n\geq 1$ such that $\ell(n)$ is finite, we have a sequence of  strong lower record values  
$$(Y_{(k)}=Y_{\ell(k)}, \ 1\leq k\leq n).$$

\Bin \textit{There are many results on probability laws of record values and record times, especially for  \textit{iid} random variables with common 
{\textit{cdf} $F$, eventually associated with the \textit{pdf} $f$ with respect to the Lebesgue measure $\lambda$}[ or \textit{iid} random variables with common mass probability functions $p$]}.\\

\Ni For now, let us focus on that case to illustrate the motivations of that study. For $n\geq 2$, the probability laws of the record values

$$
(Y^{(k)}, \ 1\leq k\leq n) \ \ and \ \  (Y_{(k)}, \ 1\leq k\leq n)
$$

\Bin are well-known. General and particular fine results on records can be found in \cite{nevzorov} and in \cite{ahsanullah_05} for example. It is of the greatest importance to remark that the probability laws of the records values heavily depend on the record time. So, the most complete approach in addressing the probability laws of records consists in giving the joint probability law of the record times and values

$$
G_n=(Y^{(1)}, \cdots, Y^{(n)}, u(1), \cdots, u(n)),
$$  

\Bin whenever $u(n)$ is finite, as given in \cite{nevzorov} (Lecture 17, Formula 17.1, page 76), and in \cite{ahsanloGL2019}, for example. That general joint law is necessary and enough to derive any measurable function $G_n$ including its margins: the joint low of record values law, the joint upper of record values, etc.\\

\Bin As far as the simultaneous joint laws of lower and upper records of $n$ records is concerned, the simultaneous joint lower-upper record values (\textit{sjlu})

$$
Z_{n}=(Y_1, \ Y_{(2)}, \cdots,  \ Y_{(n)}, Y^{(2)}, \cdots, Y^{(n)})
$$

\Bin has not been thoroughly studied. A brief study is done in \cite{arnold}, page 274, for example. But, in deriving their general formula, it is clear that the authors did no use the different possible positioning of lower record times with respect to the upper record times. The conclusion that, given the first observation 
$X_1$, the lower records and the upper records are independent is risky. This makes sense if we assume that, given $X_1=x$, we have all the upper record times
$u(1)<u(2)< \cdots < u(n)$ and next all the lower record times ($u(n)<\ell(1)<\cdots<\ell(n)$) or the reverse case where the lower record times come before the upper record times, and only in that case, the lower record values and the upper cases do not influence each other. But, these two cases are not the only ones possible as it can be seen in Figures \ref{fig1} and \ref{fig1b}. Theoretically, we have to consider all the permutations of 

$$
T_n=(\ell(1), \ \ell(2), \  \cdots, \ \ell(n), \ u(1), \ u(2), \  \cdots, \ u(n)).
$$

\Bin But any permutation in which the lower record times 

$$
(\ell(1), \ \ell(2), \  \cdots, \ell(j), \cdots, \ \ell(n))
$$

\Bin  or the upper record times 

$$
(u(1), \ u(2), \  \cdots, u(j), \cdots, \ \ u(n))
$$ 

\Bin are not given in the order in $j \in \{1,\cdots,n\}$ is an empty set. Any other permutation of the record times is possible. As a result, we will have to deal with combinatorial arguments and the final result would not be as simple as in Formula 8.4.1 in \cite{arnold}. If $n=2$, we have exactly the two simple cases described earlier and we have other cases for $n\geq 3$. In summary, a rigorous investigating of the simultaneous joint lower-upper probability law should be done through the more general probability with record times, i.e., $(Z_n, \ T_n)$, following the approach is \cite{nevzorov}.\\

\Ni The situation is not the same for the records values. Indeed, for strong records, we surely have

$$
Y_1\equiv Y_{(1)}=Y^{(1)}<Y_{(2)}<\cdots Y_{(n)} < Y^{(2)}<\cdots < Y^{(n)}.
$$

\Bin However, this ordering does not extend to the lower and upper record times as we will see below.

\Ni Here, we do not engage to find a closed-form expression of \textsl{sjlw} record values \textit{pdf} which necessarily uses combinatorial methods. Rather, we are motivated to have a complete investigations for small values of $n$, for example $n=2$ and $n=3$ for strong records for both absolutely continuous and discrete distributions and to learn from those results towards the general cases. In particular, we will learn how the  positioning of the lower record times and the upper record times in relation to each other will influence the final result.\\

\Ni We already mentioned that the two simple cases explained above are the only ones for $n=2$ records and Formula 8.4.1 in \cite{arnold} is justified. We will see that our results below, will rediscover the same results in \cite{arnold} for $2$ records. But, for $n=3$, we will see that the situation is more complicated. From the cases of $n=2$ or $n=3$ records, we will have a clear way to get the probability law $Z_n$, for any $n\geq 2$, progressively but at the cost of lengthy additive terms. So, we will fully explain the method of deriving the law of $Z_n$ for $n \in \{2,3\}$ and draw conclusions to be used in the general case.\\

\Ni To more justify our results, we will derive marginal probability laws of $Z_n$, i.e. the probability laws of $Y^{(2)}$, $Y^{(3)}$ and $(Y^{(2)}, \ Y^{(3)})$ as already found in usual record theory books and in \cite{ahsanloGL2019}.\\

\Ni We organize the rest of the paper as follows. In Section \ref{sec_ac_records}, we derive the probability law of $Z_n$ for absolutely continuous \textit{iid} random variables in two subsections, the first focusing on $n=2$ and the second on $n=3$. In Section \ref{sec_dis_records}, we do the same for discrete \textit{iid} random variables. In a last concluding section, we will draw important facts towards the general law. We recall that we deal with strong records in that paper.\\

\Ni Let us suppose that the $n$-th upper and lower record times exist and let us set for $n\geq 2$,

$$
Z_n=\left(Y_1,Y_{(2)},\cdots,Y_{(n)}, Y^{(2)},\cdots,Y^{(n)}\right). 
$$

\Bin We want to find the law of $Z_n$. In this first essay on the topic, we will study the records in an absolutely continuous frame with $f$ as the common \textit{pdf} in a section. In a second one, we focus on discrete records.

\section{Absolutely continuous records} \label{sec_ac_records}

\noindent Let us begin for the particular  case $n=2$. \\
 
\subsection{Probability law of the simultaneous joint lower-upper up to two records} \label{sec_01_ulrecords_chap2} $ $\\

\Bin Below, We state the \textit{pdf}, give the proof and derive known results as means of verification of the results and get the \textit{pdf} of $(Y_{(p)},Y^{(q)})$ for
$2\leq p,q\leq n$, when these records are defined. \\

\subsubsection{Finding the pdf of $Z_2$}

\begin{proposition} \label{prop_1} Let $Y_1$, $Y_2$, $\cdots$ be a sequence of independent real random variables defined on the same probability space $\left(\Omega, \mathcal{A}, \mathbb{P}\right)$ with common \textit{pdf} $f$ and suppose that $u(2)$ and $\ell(2)$ are finite. Then $Z_2$ has the following 
\textit{pdf}:

$$
f_{Z_2}(y_1,y_2,y_3)=f(y_1)f(y_2)f(y_3) \biggr(\frac{1}{1-F(y_1)}+ \frac{1}{F(y_1)}\biggr) 1_{(y_2<y_1<y_3)}.
$$
\end{proposition}

\Bin \textbf{Proof of Proposition \ref{prop_1}}. We simplify the writing $(t \in [y_i-dy_i/2, \ y_i+dy_i/2])$ as $t \in y_i^{\pm}$, $y_i^{+}=y_i+dy_i/2$ and $y_i^{-}=y_i-dy_i/2$, $i\geq 1$. Let $y=(y_1,y_2,y_3)\in \mathbb{R}^3$ and let us put

$$
A(dy):=\left(Y_1\in y_1^{\pm}, Y_{(2)}\in y_2^{\pm},  Y^{(2)}\in y_3^{\pm}\right),
$$

\Bin with $dy=(dy_1, dy_2, dy_3)$, \ $dy_i>0$, \ $i=1,2,3$. Put $\Delta=dy_1 \ dy_2 \ dy_3$. By definition, $Z_2$ has a \textit{pdf} if and only if\\

$$
\lim_{dy_1\rightarrow 0, dy_2\rightarrow 0, dy_3\rightarrow 0} \ \frac{\mathbb{P}(A(dy))}{\Delta}
$$

\Bin exists and hence the \textit{pdf} $f_{Z_2}$ is

$$
f_{Z_2}(y_1,y_2,y_3)=\lim_{dy_1\rightarrow 0, dy_2\rightarrow 0, dy_3\rightarrow 0} \ \frac{\mathbb{P}(A(dy))}{\Delta}, \ (y_1,y_2,y_3)\in \mathbb{R}^3.
$$

\Bin When dealing both with lower and upper records, and since an observation cannot be repeated with continuous random variables:\\

\Ni (a) the event $(u(2)=\ell(2))$ is negligible;\\

\Ni (b) the records are necessarily strong;\\

\Ni (c) $Y_2$ is either the second upper record (and $u(2)=2$) or the second lower record (and $\ell(2)=2$).\\

\Ni Now, we have 

$$
A(dy)=\biggr(A(dy) \cap (u(2)<\ell(2))\biggr) + \biggr(A(dy) \cap (\ell(2)<u(2))\biggr). 
$$

\begin{figure}
	\centering
		\includegraphics[width=.80\textwidth]{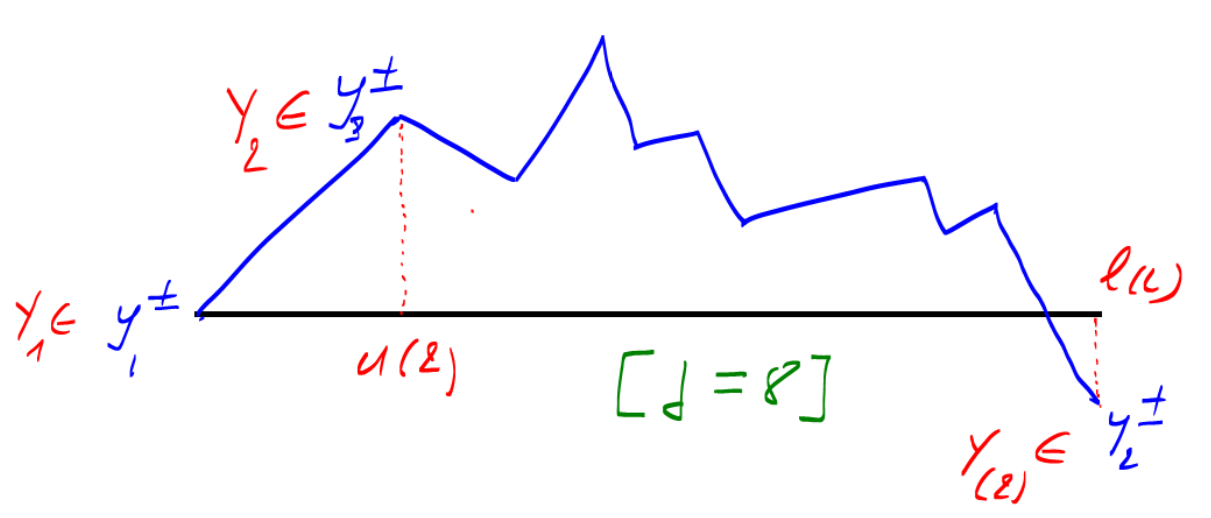}
	\caption{First lower and upper record times, where $Y_2$ is an upper record value}
	\label{fig1}
\end{figure}

\begin{figure}
	\centering
		\includegraphics[width=.80\textwidth]{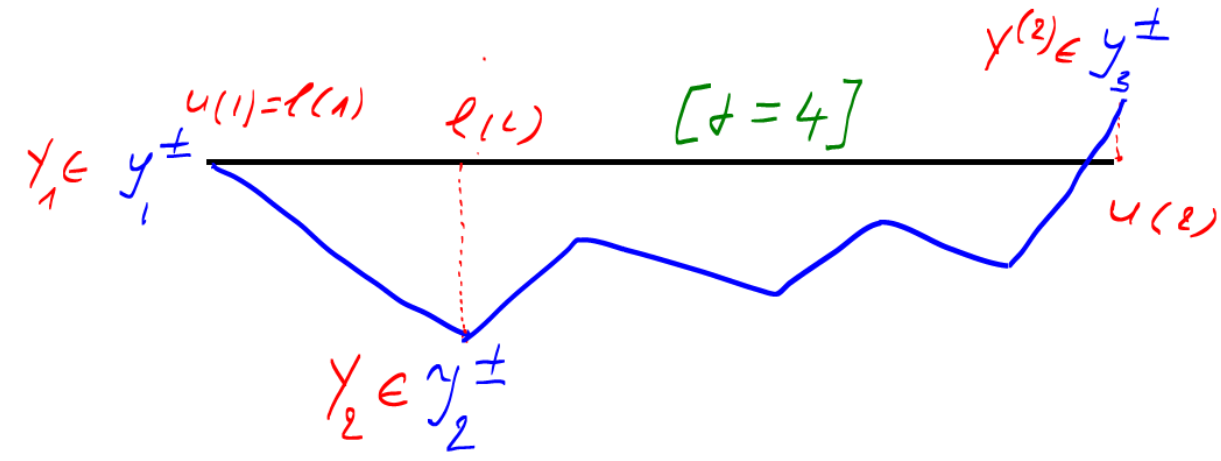}
	\caption{First lower and upper record times, where $Y_2$ is a lower record value}
	\label{fig1b}
\end{figure}

\Bin Let $N$ be the inter-record time between the two second record times. On $H_1=(A(dy) \cap (u(2)<\ell(2))$, for $dy_i$ with small norm enough, as illustrated in Fig. \ref{fig1}, we have that $y_2<y_1<y_3$ and for

$$
G_0=\biggr(Y_1 \in y_1^{\pm}, \ Y_2 \in y_3^{\pm}, \ Y_3 \in y_2^{\pm}\biggr),
$$

\Bin we have

$$
(\omega \in H_1)\cap (N=0) \Rightarrow (\omega \in G_{0}),
$$

\Bin and for $j\geq 1$,

\begin{eqnarray*}
(\omega \in H_1)\cap (N=j) &\Rightarrow& \biggr(\omega \in (Y_1 \in y_1^{\pm}) \bigcap (Y_2 \in y_3^{\pm})\cap(Y_{2+h}\geq y_1^{+}, \ 1\leq h\leq j)\\
&&\ \ \ \bigcap \ \ ( Y_{3+j} \in y_2^{\pm})\biggr).
\end{eqnarray*}

\Bin Conversely, for $y_2<y_1<y_3$, for  

\begin{eqnarray*}
G_{j}=(Y_1 \in y_1^{\pm}) \cap (Y_2 \in y_3^{\pm})\cap(Y_{2+j}\geq y_1^{+}, \ 1\leq h\leq j)\cap( Y_{3+j} \in y_2^{\pm}), \ j\geq 1,
\end{eqnarray*}

\Bin we have for any $j\geq 0$,  $\omega \in G_j$ implies that $\omega \in A(dy)\cap (N=j) \cap (u(2)<\ell(2))$ and hence

$$
G_j=H_1\cap (N=j), \ j\geq 0.
$$

\Bin But

$$
\mathbb{P}(G_j)= (F(y_1^{+})-F(y_1^{-}) ((F(y_3^{+})-F(y_3^{-})) (1-F(y_1^{+}))^j (F(y_2^{+})-F(y_2^{-})), \ j\geq 0.
$$

\Bin Finally we have

$$
\frac{\mathbb{P}(H_1)}{\Delta}=\sum_{j\geq 0} \frac{\mathbb{P}(G_j)}{\Delta}\rightarrow \frac{f(y_1)f(y_2)f(y_3)}{F(y_1)} 1_{(y_2<y_1<y_3)}. 
$$

\Bin We treat $(\omega \in H_2)\cap (N=j)$ in the same manner (see  Fig. \ref{fig1b}) where $H_2=(A(dy) \cap (u(2)>\ell(2))$ and the $G_j$'s are replaced by, for $y_2<y_1<y_3$,

$$
G_0=\biggr(Y_1 \in y_1^{\pm}, \ Y_2 \in y_2^{\pm}, \ Y_3 \in y_3^{\pm}\biggr),
$$

\Bin and

$$
G_{j}=(Y_1 \in y_1^{\pm}) \cap (Y_2 \in y_2^{\pm})\cap(Y_{2+j}\leq y_1^{-}, \ 1\leq h\leq j)\cap( Y_{3+j} \in y_3^{\pm}), \ j\geq 1,
$$

\Bin and 

$$
\mathbb{P}(G_j)= (F(y_1^{+})-F(y_1^{-})) (F(y_2^{+})-F(y_2^{-})) (F(y_1^{-})^j (F(y_3^{+})-F(y_3^{-})), \ j\geq 0.
$$

\Bin We conclude that

$$
\frac{\mathbb{P}(H_2)}{\Delta}=\sum_{j\geq 0} \frac{\mathbb{P}(G_j)}{\Delta}\rightarrow \frac{f(y_1)f(y_2)f(y_3)}{1-F(y_1)} 1_{(y_2<y_1<y_3)}. 
$$

\Bin The proof is over. $\blacksquare$\\

\subsubsection{Derivation of known results and of the second record values (lower and upper)} $ $ \\

\noindent \textbf{(i) Let us rediscover the law of $\left(Y_1,Y^{(2)}\right)$ whose \textit{pdf} is}:

\begin{eqnarray*}
f_{(Y_1,Y^{(2)})}(y_1,y_3) &=& \int f_{Z_2}(y_1,y_2,y_3)dy_2 \\
&=& 1_{(y_1<y_3)} \left(\frac{1}{F(y_1)} + \frac{1}{1-F(y_1)}\right) f(y_1)f(y_3) \int_{-\infty}^{y_1}\ f(y_2) dy_2 \\
&=& 1_{(y_1<y_3)} \left(\frac{1}{F(y_1)} + \frac{1}{1-F(y_1)}\right) f(y_1)f(y_3)F(y_1) \\
&=& 1_{(y_1<y_3)} \left(1 + \frac{F(y_1)}{1-F(y_1)}\right) f(y_1)f(y_3) \\ 
&=& \frac{f(y_1)}{1-F(y_1)} f(y_3)1_{(y_1<y_3)}.
\end{eqnarray*}

\Bin So, if $r(x)=f(x)/(1-F(x))$, $x\in ]\textit{\text{lep}}(F), \ \textit{\text{uep}}(F)[$, then

$$
f_{(Y_1,Y^{(2)})}(x,y) = r(x) f(y)1_{(x<y)}
$$

\Bin and this confirms a known result. \\

\noindent \textbf{(ii) We can derive the law of $\left(Y_{(2)},Y^{(2)}\right)$ as given below}:

\begin{corollary} \label{low_upp_ac_01}

\Ni The \textit{pdf} of $\left(Y_{(2)},Y^{(2)}\right)$ is then 

$$
f_{(Y_{(2)},Y^{(2)})}(y,z) = f(y)f(z)\left(\log \frac{F(z)}{F(y)} + \log \frac{1-F(y)}{1-F(z)}\right)1_{(y<z)}. 
$$

\end{corollary}

\Bin \textbf{Proof of Corollary \ref{low_upp_ac_01}}. We have

\begin{eqnarray*}
f_{(Y_{(2)},Y^{(2)})}(y_2,y_3) &=& \int f_{Z_2}(y_1,y_2,y_3)dy_1 \\
&=& 1_{(y_2<y_3)}f(y_2)f(y_3) \int_{y_2}^{y_3} \left(\frac{f(y_1)}{F(y_1)} + \frac{f(y_1)}{1-F(y_1)}\right) dy_1 \\
&=& f(y_2)f(y_3)\left((R_1(y_3)-R_1(y_2)) + (R_2(y_3)-R_2(y_2))\right)1_{(y_2<y_3)},
\end{eqnarray*}

\Bin with $R_1(t)=\log F(t)$ \ and \ $R_2(t)=-\log(1-F(t))$. So we finally have

$$
f_{(Y_{(2)},Y^{(2)})}(y,z) = f(y)f(z)\left(\log \frac{F(z)}{F(y)} + \log \frac{1-F(y)}{1-F(z)}\right)1_{(y<z)}
$$

\Bin and the proof is over. $\blacksquare$ \\

\subsection{Probability law of the simultaneous joint lower-upper up to three records} \label{sec_02_ulrecords_chap2}  $ $\\

\subsubsection{Finding the \textit{pdf} of $Z_3$}

\begin{proposition} \label{prop_2}
Let $Y_1$, $Y_2$, $\cdots$ be a sequence of independent real random variables defined on the same probability space $\left(\Omega, \mathcal{A}, \mathbb{P}\right)$ with common \textit{pdf} $f$ and suppose that $u(3)$ and $\ell(3)$ are finite. Then $Z_3$ has the following \textit{pdf}:

$$
f_{Z_3}(x,y_1,y_2,z_1,z_2)=f(x)f(y_1)f(y_2)f(z_1)f(z_2) L(x,y_1,z_1) 1_{(y_2<y_1<x<z_1<z_2)},
$$

\Bin for

\begin{eqnarray*}
&&L(x,y_1,z_1)\\
&&=\frac{1}{(1-[F(z_1)-F(x)])F(x)F(y_1)}\\
&&+\frac{1}{(1-[F(z_1)-F(x)])(1-[F(z_1)-F(y_1)])(1-F(z_1))}\\
&&+\frac{1}{(1-[F(z_1)-F(x)])(1-[F(z_1)-F(y_1)])F(y_1)}\\
&&+\frac{1}{(1-[F(x)-F(y_1)])(1-[F(z_1)-F(y_1)])(1-F(z_1))}\\
&&+\frac{1}{(1-[F(x)-F(y_1)])(1-[F(z_1)-F(y_1)])F(y_1)}\\
&&+\frac{1}{(1-[F(x)-F(y_1)])(1-F(x))(1-F(z_1))}.
\end{eqnarray*}

\end{proposition}

\Bin \textbf{Proof of Proposition \ref{prop_2}}. We use similar methods to those in Section \ref{sec_01_ulrecords_chap2} but the situation a little more complex. We define

$$
A(d):=\left(Y_1\in x^{\pm}, Y_{(2)}\in y_1^{\pm}, Y_{(3)}\in y_2^{\pm}, Y^{(2)}\in z_1^{\pm}, Y^{(3)}\in z_2^{\pm}\right),
$$

\Bin with $d=(dx,dy_1,dy_2,dz_1,dz_2)$ and $\Delta=dx \times dy_1 \times dy_2 \times dz_1 \times dz_2$. We have to consider all the 24 orderings of
$(u(2),u(3),\ell(2),\ell(3))$. Fortunately $u(3)$ (resp. $\ell(3)$) cannot come before $u(2)$ (resp. $\ell(2)$). It will remains six orderings

\begin{eqnarray*}
O_1&=&(u(2)<u(3)<\ell(2)<\ell(3))\\
O_2&=&(u(2)<\ell(2)<\ell(3)<u(3))\\
O_3&=&(u(2)<\ell(2)<u(3)<\ell(3))\\
O_4&=&(\ell(2)<u(2)<\ell(3)<u(3))\\
O_5&=&(\ell(2)<u(2)<u(3)<\ell(3))\\
O_6&=&(\ell(2)<\ell(3)<u(2)<u(3)).
\end{eqnarray*} 

\Bin Let $N(h)$, $h \in \{1,2,3\}$ the three inter-record times. We have to decompose each $H_i=A(d)\cap O_i$, $i \in \{1,\cdots,6\}$ into

$$
H_i=\sum_{i_1\geq 0, i_2\geq 0, i_3\geq 0} H_i \cap (N(1)=i_1,N(2)=i_2,N(3)=i_3). 
$$

\begin{figure}
	\centering
		\includegraphics[width=1.00\textwidth]{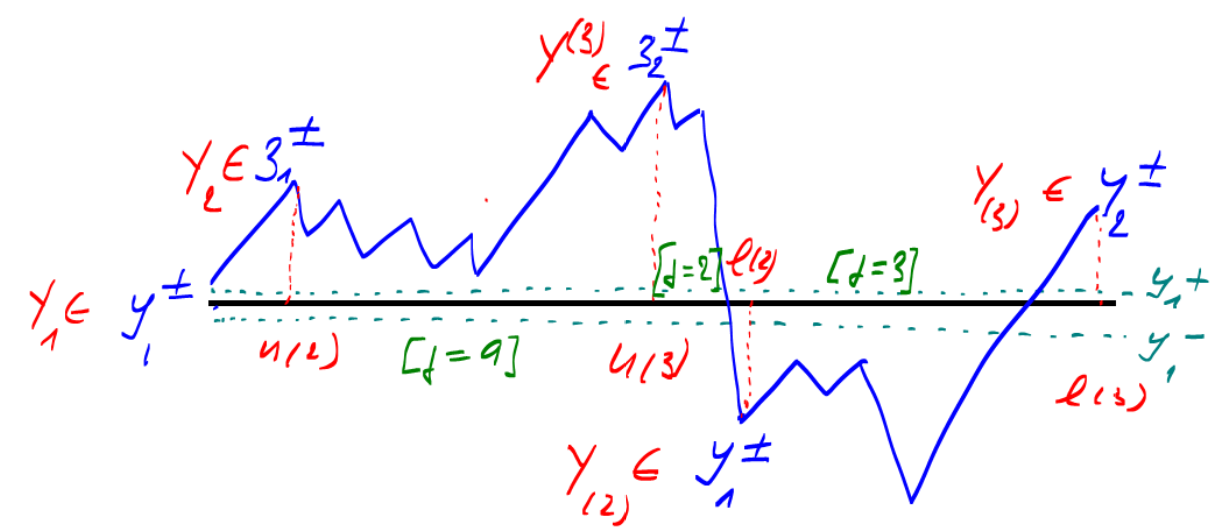}
	\caption{Case $u(2)<u(3)<\ell(2)<\ell(3)$}
	\label{fig2}
\end{figure}

\Bin As above, we are going to describe $G_{i,i_1,i_2,i_3}=H_i \cap (N(1)=i_1,N(2)=i_2,N(3)=i_3)$. We have six cases to deal with. So, we fully explain one cases and let the reader check the other five cases. From Fig. \ref{fig2}, we have the following facts (for $d$ small enough in norm):\\

\Ni (a) For $h \in ]u(2),u(3)[$, $Y_h>x^{+}$ otherwise it would be the second lower record value, and $Y_h<z_1^{-}$ otherwise it would be the third upper record;\\

\Ni (b) For $h \in ]u(3),\ell(2)[$, $Y_h>x^{+}$ otherwise it would be the second lower record (but it may exceed $z_2^{+}$ without any consequence);\\

\Ni (c) For $h \in ]\ell(2),\ell(3)[$, $Y_h>y_1^{+}$ otherwise it would be the third lower record (but it is not bounded above).\\

\Ni Hence, we get 
$$
\frac{\mathbb{P}(G_{1,i_1,i_2,i_3})}{\Delta}\rightarrow f(x)f(y_1)f(y_2)f(z_1)f(z_2) (F(z_1)-F(x))^{i_1} (1-F(x))^{i_2} (1-F(y_1))^{i_3},
$$

\Bin and hence

$$
\frac{\mathbb{P}(H_1)}{\Delta}\rightarrow \frac{f(x)f(y_1)f(y_2)f(z_1)f(z_2)}{(1-[F(z_1)-F(x)])F(x)F(y_1)}.
$$

\Bin For the five other cases, we have

$$
\frac{\mathbb{P}(G_{2,i_1,i_2,i_3})}{\Delta}\rightarrow f(x)f(y_1)f(y_2)f(z_1)f(z_2) (F(z_1)-F(x))^{i_1} (F(z_1)-F(y_1))^{i_2} F(z_1)^{i_3},
$$

$$
\frac{\mathbb{P}(G_{3,i_1,i_2,i_3})}{\Delta}\rightarrow f(x)f(y_1)f(y_2)f(z_1)f(z_2) (F(z_1)-F(x))^{i_1} (F(z_1)-F(y_1))^{i_2} (1-F(y_1))^{i_3},
$$

$$
\frac{\mathbb{P}(G_{4,i_1,i_2,i_3})}{\Delta}\rightarrow f(x)f(y_1)f(y_2)f(z_1)f(z_2) (F(x)-F(y_1))^{i_1} (F(z_1)-F(y_1))^{i_2} F(z_1)^{i_3},
$$

$$
\frac{\mathbb{P}(G_{5,i_1,i_2,i_3})}{\Delta}\rightarrow f(x)f(y_1)f(y_2)f(z_1)f(z_2) (F(x)-F(y_1))^{i_1} (F(z_1)-F(y_1))^{i_2} (1-F(y_1))^{i_3},
$$

$$
\frac{\mathbb{P}(G_{6,i_1,i_2,i_3})}{\Delta}\rightarrow f(x)f(y_1)f(y_2)f(z_1)f(z_2) (F(x)-F(y_1))^{i_1} F(x)^{i_2} F(z_1)^{i_3},
$$

\Ni These facts when put together lead to

$$
f_{Z_3}(x,y_1,y_2,z_1,z_2)= f(x)f(y_1)f(y_2)f(z_1)f(z_2) L(x,y_1,z_1)1_{(y_2<y_1<x<z_1<z_2)}
$$

\Bin with

\begin{eqnarray*}
&&L(x,y_1,z_1)\\
&&=\frac{1}{(1-[F(z_1)-F(x)])F(x)F(y_1)}\\
&&+\frac{1}{(1-[F(z_1)-F(x)])(1-[F(z_1)-F(y_1)])(1-F(z_1))}\\
&&+\frac{1}{(1-[F(z_1)-F(x)])(1-[F(z_1)-F(y_1)])F(y_1)}\\
&&+\frac{1}{(1-[F(x)-F(y_1)])(1-[F(z_1)-F(y_1)])(1-F(z_1))}\\
&&+\frac{1}{(1-[F(x)-F(y_1)])(1-[F(z_1)-F(y_1)])F(y_1)}\\
&&+\frac{1}{(1-[F(x)-F(y_1)])(1-F(x))(1-F(z_1))}.
\end{eqnarray*}

\Bin

\subsubsection{Derivation of known results and of the third record values (lower and upper)} $ $\\

\noindent \textbf{(i) Checking of results by comparing with the known \textit{pdf} of $\left(Y_1, Y^{(2)}, Y^{(3)}\right)$}.\\

\Ni We have 

\begin{eqnarray*}
f_{(Y_1,Y^{(2)},Y^{(3)})}(x,z_1,z_2) &=& \int \int f_{Z_3}(x,y_1,y_2,z_1,z_2)dy_1 dy_2 \\
&=:& 1_{(x<z_1<z_2)} f(x)f(z_1)f(z_2)\left(I_1 + I_2 + I_3 + I_4 + I_5 + I_6\right).
\end{eqnarray*}

\Ni Let us denote, for $(x,z_1,z_2)$ fixed, by $D_x := D_{(x,z_1,z_2)} = (y_2<y_1<x)$ as the section of the domain of $Z_3$ at $(x,z_1,z_2)$. So

\begin{eqnarray*}
I_1 &=& \int \int_{D_x} \frac{f(y_1)f(y_2)}{(1-[F(z_1)-F(x)])F(x)F(y_1)}dy_1 dy_2 \\
&=& \frac{1}{(1-[F(z_1)-F(x)])F(x)} \int \int_{D_x} \frac{f(y_1)f(y_2)}{F(y_1)}dy_1 dy_2  \\
&=& \frac{1}{(1-[F(z_1)-F(x)])F(x)} \int_{-\infty}^{x} \frac{f(y_1)}{F(y_1)} \left(\int_{-\infty}^{y_1} f(y_2) dy_2\right) dy_1 \\
&=& \frac{1}{(1-[F(z_1)-F(x)])}. 
\end{eqnarray*}

\begin{eqnarray*}
I_2 &=& \int \int_{D_x} \frac{f(y_1)f(y_2)}{(1-[F(z_1)-F(x)])(1-[F(z_1)-F(y_1)])(1-F(z_1))}dy_1 dy_2 \\
&=& \frac{1}{(1-[F(z_1)-F(x)])(1-F(z_1))} \int \int_{D_x} \frac{f(y_1)f(y_2)}{1-[F(z_1)-F(y_1)]}dy_1 dy_2  \\
&=& \frac{1}{(1-[F(z_1)-F(x)])(1-F(z_1))} \int_{-\infty}^{x} \frac{f(y_1)}{1-[F(z_1)-F(y_1)]} \left(\int_{-\infty}^{y_1} f(y_2) dy_2\right) dy_1 \\
&=& \frac{1}{(1-[F(z_1)-F(x)])(1-F(z_1))} \int_{-\infty}^{x} \frac{f(y_1) F(y_1)}{1-[F(z_1)-F(y_1)]} dy_1 \\ 
&=& \frac{1}{(1-[F(z_1)-F(x)])(1-F(z_1))} \int_{-\infty}^{x} \frac{f(y_1)\left\{1-[F(z_1)-F(y_1)]\right\} - f(y_1)(1-F(z_1))}{1-[F(z_1)-F(y_1)]} dy_1 \\
&=& \frac{1}{(1-[F(z_1)-F(x)])(1-F(z_1))} \left\{\int_{-\infty}^{x} f(y_1) dy_1 - \int_{-\infty}^{x} \frac{f(y_1)(1-F(z_1))}{1-[F(z_1)-F(y_1)]} dy_1\right\} \\
&=& \frac{1}{(1-[F(z_1)-F(x)])(1-F(z_1))} \left\{F(x) - (1-F(z_1))\log \frac{1-[F(z_1)-F(x)]}{1-F(z_1)} \right\} \\
&=& \frac{F(x)}{(1-[F(z_1)-F(x)])(1-F(z_1))} - \frac{\log \left(1-[F(z_1)-F(x)]\right) - \log \left(1-F(z_1)\right)}{1-[F(z_1)-F(x)]}.
\end{eqnarray*}

\begin{eqnarray*}
I_3 &=& \int \int_{D_x} \frac{f(y_1)f(y_2)}{(1-[F(z_1)-F(x)])(1-[F(z_1)-F(y_1)])F(y_1)}dy_1 dy_2 \\
&=& \frac{1}{1-[F(z_1)-F(x)]} \int \int_{D_x} \frac{f(y_1)f(y_2)}{(1-[F(z_1)-F(y_1)])F(y_1)}dy_1 dy_2  \\
&=& \frac{1}{1-[F(z_1)-F(x)]} \int_{-\infty}^{x} \frac{f(y_1)}{(1-[F(z_1)-F(y_1)])F(y_1)} \left(\int_{-\infty}^{y_1} f(y_2) dy_2\right) dy_1 \\
&=& \frac{1}{1-[F(z_1)-F(x)]} \int_{-\infty}^{x} \frac{f(y_1)}{1-[F(z_1)-F(y_1)]} dy_1 \\ 
&=& \frac{\log \left(1-[F(z_1)-F(x)]\right) - \log \left(1-F(z_1)\right)}{1-[F(z_1)-F(x)]}.
\end{eqnarray*}

\begin{eqnarray*}
I_4 &=& \int \int_{D_x} \frac{f(y_1)f(y_2)}{(1-[F(x)-F(y_1)])(1-[F(z_1)-F(y_1)])(1-F(z_1))}dy_1 dy_2 \\
&=& \frac{1}{(1-F(z_1))} \int \int_{D_x} \frac{f(y_1)f(y_2)}{(1-[F(x)-F(y_1)])(1-[F(z_1)-F(y_1)])}dy_1 dy_2  \\
&=& \frac{1}{(1-F(z_1))} \int_{-\infty}^{x} \frac{f(y_1)}{(1-[F(x)-F(y_1)])(1-[F(z_1)-F(y_1)])} \left(\int_{-\infty}^{y_1} f(y_2) dy_2\right) dy_1 \\
&=& \frac{1}{(1-F(z_1))} \int_{-\infty}^{x} \frac{f(y_1) F(y_1)}{(1-[F(x)-F(y_1)])(1-[F(z_1)-F(y_1)])} dy_1. \\ 
\end{eqnarray*}

\Bin Hence we use the following decomposition 

\begin{eqnarray*}
\frac{F(y_1)}{(1-[F(x)-F(y_1)])(1-[F(z_1)-F(y_1)])} &=& \frac{1-F(z_1)}{F(x)-F(z_1)}\times \frac{1}{1-[F(z_1)-F(y_1)]} \\ 
&-& \frac{1-F(x)}{F(x)-F(z_1)}\times \frac{1}{1-[F(x)-F(y_1)]}
\end{eqnarray*}

\Bin and so

\begin{eqnarray*}
I_4 &=& \frac{1}{(1-F(z_1))} \biggr(\frac{1-F(z_1)}{F(x)-F(z_1)} \int_{-\infty}^{x} \frac{f(y_1)}{1-[F(z_1)-F(y_1)]} dy_1 \\
&-& \frac{1-F(x)}{F(x)-F(z_1)}\int_{-\infty}^{x}\frac{f(y_1)}{1-[F(x)-F(y_1)]} dy_1 \biggr) \\
&=& \frac{1}{(1-F(z_1))} \biggr(\frac{1-F(z_1)}{F(x)-F(z_1)} \left(\log \left(1-[F(z_1)-F(x)]\right) - \log \left(1-F(z_1)\right)\right)  \\
&+& \frac{1-F(x)}{F(x)-F(z_1)} \log \left(1-F(x)\right) \biggr) \\
&=& \frac{1}{(1-F(z_1))(F(x)-F(z_1))} \biggr( (1-F(z_1))\biggr(\log \left(1-[F(z_1)-F(x)]\right) - \log \left(1-F(z_1)\right)\biggr) \\ 
&+& (1-F(x))\log \left(1-F(x)\right)\biggr).
\end{eqnarray*}

\begin{eqnarray*}
I_5 &=& \int \int_{D_x} \frac{f(y_1)f(y_2)}{(1-[F(x)-F(y_1)])(1-[F(z_1)-F(y_1)])F(y_1)}dy_1 dy_2 \\
&=& \int_{-\infty}^{x} \frac{f(y_1)}{(1-[F(x)-F(y_1)])(1-[F(z_1)-F(y_1)])F(y_1)}\left(\int_{-\infty}^{y_1} f(y_2) dy_2\right)dy_1  \\
&=& \int_{-\infty}^{x} \frac{f(y_1)}{(1-[F(x)-F(y_1)])(1-[F(z_1)-F(y_1)])}dy_1  \\
&=& \frac{1}{F(x)-F(z_1)} \int_{-\infty}^{x} \frac{f(y_1)(F(x)-F(z_1))}{(1-[F(x)-F(y_1)])(1-[F(z_1)-F(y_1)])} dy_1 \\
&=& \frac{1}{F(x)-F(z_1)} \biggr(\int_{-\infty}^{x} \frac{f(y_1)}{1-[F(x)-F(y_1)]} dy_1 - \int_{-\infty}^{x} \frac{f(y_1)}{1-[F(z_1)-F(y_1)]} dy_1\biggr)  \ (L4)\\
&=& \frac{\log \left(1-F(z_1)\right) - \log \left(1-F(x)\right) - \log \left(1-[F(z_1)-F(x)]\right)}{F(x)-F(z_1)} 
\end{eqnarray*}

\Bin where we used in Line (L4) the simple fact

$$
\frac{(F(x)-F(z_1))}{(1-[F(x)-F(y_1)])(1-[F(z_1)-F(y_1)])} = \frac{1}{1-[F(x)-F(y_1)]} - \frac{1}{1-[F(z_1)-F(y_1)]}
$$

\begin{eqnarray*}
I_6 &=& \int \int_{D_x} \frac{f(y_1)f(y_2)}{(1-[F(x)-F(y_1)])(1-F(x))(1-F(z_1))}dy_1 dy_2 \\
&=& \frac{1}{(1-F(x))(1-F(z_1))} \int \int_{D_x} \frac{f(y_1)f(y_2)}{1-[F(x)-F(y_1)]}dy_1 dy_2  \\
&=& \frac{1}{(1-F(x))(1-F(z_1))} \int_{-\infty}^{x} \frac{f(y_1)}{1-[F(x)-F(y_1)]} \left(\int_{-\infty}^{y_1} f(y_2) dy_2\right) dy_1 \\
&=& \frac{1}{(1-F(x))(1-F(z_1))} \int_{-\infty}^{x} \frac{f(y_1)F(y_1)}{1-[F(x)-F(y_1)]} dy_1 \\
&=& \frac{1}{(1-F(x))(1-F(z_1))} \biggr(\int_{-\infty}^{x} f(y_1) dy_1 \\
&-& \int_{-\infty}^{x} (1-F(x))\frac{f(y_1)}{1-[F(x)-F(y_1)]} dy_1 \biggr) \\
&=& \frac{1}{(1-F(x))(1-F(z_1))} \biggr( F(x) + (1-F(x))\log(1-F(x)) \biggr) \\
&=& \frac{F(x)}{(1-F(x))(1-F(z_1))} + \frac{\log(1-F(x))}{1-F(z_1)}
\end{eqnarray*}

\Bin Hence by doing the simple computations, we arrive at 

$$
I_1 + I_2 + I_3 + I_4 + I_5 + I_6 = \frac{1}{(1-F(x))(1-F(z_1))}
$$

\Bin and hence

\begin{eqnarray*}
f_{(Y_1,Y^{(2)},Y^{(3)})}(x,z_1,z_2) &=&  \frac{f(x)f(z_1)f(z_2)}{(1-F(x))(1-F(z_1))} 1_{(x<z_1<z_2)} \\
&=& r(x) r(z_1) f(z_2) 1_{(x<z_1<z_2)} 
\end{eqnarray*}

\Bin and this confirms a known result, with

$$
r(t)=f(t)/(1-F(t)), \ \  t\in ]\textit{\text{lep}}(F), \ \textit{\text{uep}}(F)[.
$$ 

\Bin  \textbf{(ii) Let us derive the laws of $\left(Y_{(p)},Y^{(q)}\right)$, $2\leq p,q \leq 3$}.\\

\Ni \textbf{Let us begin by checking the law of \textbf{$\left(Y_{(2)},Y^{(2)}\right)$}, which has already been found in Corollary \ref{low_upp_ac_01}.} \\

\Ni We have

\begin{eqnarray*}
f_{(Y_{(2)},Y^{(2)})}(y_1,z_1) &=& \int \int \int f_{Z_3}(x,y_1,y_2,z_1,z_2) dx dy_2 dz_2 \\
&=:& 1_{(y_1<z_1)} f(y_1)f(z_1)\left(J_1 + J_2 + J_3 + J_4 + J_5 + J_6\right).
\end{eqnarray*}

\Bin Let us denote by, for $(y_1,z_1)$ fixed,  $D_{(y_1,z_1)} = (y_2<y_1<x<z_1<z_2) =: D_x$ as the section of the domain of $Z_3$ at $(y_1,z_1)$ . So

\begin{eqnarray*}
J_1 &=& \int \int \int_{D_x} \frac{f(x)f(y_2)f(z_2)}{(1-[F(z_1)-F(x)])F(x)F(y_1)} dx dy_2 dz_2 \\
&=& \frac{1}{F(y_1)} \int_{-\infty}^{y_1} f(y_2)dy_2 \int_{z_1}^{+\infty} f(z_2)dz_2 \int_{y_1}^{z_1} \frac{f(x)}{(1-[F(z_1)-F(x)])F(x)} dx  \\
&=& \int_{y_1}^{z_1} \frac{f(x)(1-F(z_1))}{(1-[F(z_1)-F(x)])F(x)} dx.  \\
\end{eqnarray*}

\Bin But remark that 

$$
\frac{1-F(z_1)}{(1-[F(z_1)-F(x)])F(x)} = \frac{1}{F(x)} - \frac{1}{1-[F(z_1)-F(x)]}
$$

\Bin and hence

\begin{eqnarray*}
J_1 &=& \int_{y_1}^{z_1} \biggr(\frac{f(x)}{F(x)} - \frac{f(x)}{1-[F(z_1)-F(x)]}\biggr) dx  \\
&=& \log \frac{F(z_1)}{F(y_1)} + \log\left(1-[F(z_1)-F(y_1)]\right).
\end{eqnarray*}

\begin{eqnarray*}
J_2 &=& \int \int \int_{D_x} \frac{f(x)f(y_2)f(z_2)}{(1-[F(z_1)-F(x)])(1-[F(z_1)-F(y_1)])(1-F(z_1))} dx dy_2 dz_2 \\
&=& \frac{1}{(1-[F(z_1)-F(y_1)])(1-F(z_1))} \int_{-\infty}^{y_1} f(y_2)dy_2 \\
&&\int_{z_1}^{+\infty} f(z_2)dz_2 \int_{y_1}^{z_1} \frac{f(x)}{1-[F(z_1)-F(x)]} dx  \\
&=& -F(y_1)\times \frac{\log\left(1-[F(z_1)-F(y_1)]\right)}{1-[F(z_1)-F(y_1)]}.
\end{eqnarray*}

\begin{eqnarray*}
J_3 &=& \int \int \int_{D_x} \frac{f(x)f(y_2)f(z_2)}{(1-[F(z_1)-F(x)])(1-[F(z_1)-F(y_1)])F(y_1)} dx dy_2 dz_2 \\
&=& \frac{1}{(1-[F(z_1)-F(y_1)])F(y_1)} \int_{-\infty}^{y_1} f(y_2)dy_2\\
&& \int_{z_1}^{+\infty} f(z_2)dz_2 \int_{y_1}^{z_1} \frac{f(x)}{1-[F(z_1)-F(x)]} dx  \\
&=& -(1-F(z_1))\times \frac{\log\left(1-[F(z_1)-F(y_1)]\right)}{1-[F(z_1)-F(y_1)]}.
\end{eqnarray*}

\begin{eqnarray*}
J_4 &=& \int \int \int_{D_x} \frac{f(x)f(y_2)f(z_2)}{(1-[F(x)-F(y_1)])(1-[F(z_1)-F(y_1)])(1-F(z_1))} dx dy_2 dz_2 \\
&=& \frac{1}{(1-[F(z_1)-F(y_1)])(1-F(z_1))} \int_{-\infty}^{y_1} f(y_2)dy_2 \int_{z_1}^{+\infty} f(z_2)dz_2 \\
&&\int_{y_1}^{z_1} \frac{f(x)}{1-[F(x)-F(y_1)]} dx  \\
&=& -F(y_1)\times \frac{\log\left(1-[F(z_1)-F(y_1)]\right)}{1-[F(z_1)-F(y_1)]}.
\end{eqnarray*}

\begin{eqnarray*}
J_5 &=& \int \int \int_{D_x} \frac{f(x)f(y_2)f(z_2)}{(1-[F(x)-F(y_1)])(1-[F(z_1)-F(y_1)])F(y_1)} dx dy_2 dz_2 \\
&=& \frac{1}{(1-[F(z_1)-F(y_1)])F(y_1)} \int_{-\infty}^{y_1} f(y_2)dy_2 \\
&&\int_{z_1}^{+\infty} f(z_2)dz_2 \int_{y_1}^{z_1} \frac{f(x)}{1-[F(x)-F(y_1)]} dx  \\
&=& -(1-F(z_1))\times \frac{\log\left(1-[F(z_1)-F(y_1)]\right)}{1-[F(z_1)-F(y_1)]}.
\end{eqnarray*}

\begin{eqnarray*}
J_6 &=& \int \int \int_{D_x} \frac{f(x)f(y_2)f(z_2)}{(1-[F(x)-F(y_1)])(1-F(x))(1-F(z_1))} dx dy_2 dz_2 \\
&=& \frac{1}{(1-F(z_1))} \int_{-\infty}^{y_1} f(y_2)dy_2 \int_{z_1}^{+\infty} f(z_2)dz_2 \\
&&\int_{y_1}^{z_1} \frac{f(x)}{(1-[F(x)-F(y_1)])(1-F(x))} dx  \\
&=& \int_{y_1}^{z_1} \frac{f(x) F(y_1)}{(1-[F(x)-F(y_1)])(1-F(x))} dx.  \\
\end{eqnarray*}

\Bin But remark that 

$$
\frac{F(y_1)}{(1-[F(x)-F(y_1)])(1-F(x))} = \frac{1}{1-F(x)} - \frac{1}{1-[F(x)-F(y_1)]}
$$

\Bin and hence

\begin{eqnarray*}
J_6 &=& \int_{y_1}^{z_1} \biggr(\frac{f(x)}{1-F(x)} - \frac{f(x)}{1-[F(x)-F(y_1)]}\biggr) dx  \\
&=& \log \frac{1-F(y_1)}{1-F(z_1)} + \log\left(1-[F(z_1)-F(y_1)]\right).
\end{eqnarray*}

\Bin Hence by regrouping all terms, we will have

\begin{eqnarray*}
f_{(Y_{(2)},Y^{(2)})}(y_1,z_1) = f(y_1)f(z_1)\left(\log \frac{F(z_1)}{F(y_1)} + \log \frac{1-F(y_1)}{1-F(z_1)}\right) 1_{(y_1<z_1)}
\end{eqnarray*}

\Bin and that confirms the result found in Corollary \ref{low_upp_ac_01}. \\

\Ni \textbf{Now, we are going to derive the laws of $\left(Y_{(2)},Y^{(3)}\right)$, $\left(Y_{(3)},Y^{(2)}\right)$ and $\left(Y_{(3)},Y^{(3)}\right)$ in the next Corollaries.} \\

\begin{corollary} \label{low_upp_ac_02}
The \textit{pdf} of $\left(Y_{(2)},Y^{(3)}\right)$ is then given by

$$
f_{(Y_{(2)},Y^{(3)})}(y,z) = -  f(y)f(z)\biggr(\int_{y}^{z} \frac{f(x)}{F(x)(1-F(x))}\log \frac{1-F(z)}{1-F(x)} dx\biggr)1_{(y<z)}.
$$

\end{corollary}

\Bin \textbf{Proof of Corollary \ref{low_upp_ac_02}}.
We have 

\begin{eqnarray*}
f_{(Y_{(2)},Y^{(3)})}(y_1,z_2) &=& \int \int \int f_{Z_3}(x,y_1,y_2,z_1,z_2) dx dy_2 dz_1 \\
&=:& 1_{(y_1<z_2)} f(y_1)f(z_2)\left(K_1 + K_2 + K_3 + K_4 + K_5 + K_6\right).
\end{eqnarray*}

\Bin Let us denote, for $(y_1,z_2)$ fixed, by  $D_{(y_1,z_2)} = (y_2<y_1<x<z_1<z_2) =: D_x$ as the section of the domain of $Z_3$ at $(y_1,z_2)$. So

\begin{eqnarray*}
K_1 &=& \int \int \int_{D_x} \frac{f(x)f(y_2)f(z_1)}{(1-[F(z_1)-F(x)])F(x)F(y_1)} dx dy_2 dz_1 \\
&=&\frac{1}{F(y_1)} \int_{-\infty}^{y_1} f(y_2) dy_2 \int_{y_1}^{z_2} \frac{f(x)}{F(x)}\biggr(\int_{x}^{z_2} \frac{f(z_1)}{(1-[F(z_1)-F(x)])} dz_1 \biggr) dx  \\
&=& -\int_{y_1}^{z_2} \frac{f(x)}{F(x)} \log\left(1-[F(z_2)-F(x)]\right) dx.  
\end{eqnarray*}

\begin{eqnarray*}
K_2 &=& \int \int \int_{D_x} \frac{f(x)f(y_2)f(z_1)}{(1-[F(z_1)-F(x)])(1-[F(z_1)-F(y_1)])(1-F(z_1))} dx dy_2 dz_1 \\
&=& \int_{-\infty}^{y_1} f(y_2) dy_2 \int_{y_1}^{z_2} f(x)\\
&&\biggr(\int_{x}^{z_2} \frac{f(z_1)}{(1-[F(z_1)-F(x)])(1-[F(z_1)-F(y_1)])(1-F(z_1))} dz_1 \biggr) dx  \\
&=& F(y_1)\int_{y_1}^{z_2} f(x)\\
&&\biggr(\int_{x}^{z_2} \frac{f(z_1)}{(1-[F(z_1)-F(x)])(1-[F(z_1)-F(y_1)])(1-F(z_1))} dz_1 \biggr) dx.   
\end{eqnarray*}

\Bin Now we use the following relation 

\begin{eqnarray*}
&&\frac{1}{(1-[F(z_1)-F(x)])(1-[F(z_1)-F(y_1)])(1-F(z_1))} \\
&&=\frac{1}{F(x)(F(y_1)-F(x))}\biggr(\frac{1}{1-F(z_1)} - \frac{1}{1-[F(z_1)-F(x)]}\biggr) \\
&&-\frac{1}{F(y_1)(F(y_1)-F(x))}\biggr(\frac{1}{1-F(z_1)} - \frac{1}{1-[F(z_1)-F(y_1)]}\biggr).
\end{eqnarray*}

\Bin to get

\begin{eqnarray*}
K_2 &=& F(y_1)\int_{y_1}^{z_2} \frac{f(x)}{F(x)(F(y_1)-F(x))}\biggr(\left\{\frac{f(z_1)}{1-F(z_1)} - \frac{f(z_1)}{1-[F(z_1)-F(x)]}\right\}dz_1\biggr) dx \\
&-& F(y_1)\int_{y_1}^{z_2} \frac{f(x)}{F(y_1)(F(y_1)-F(x))}\biggr(\left\{\frac{f(z_1)}{1-F(z_1)} - \frac{f(z_1)}{1-[F(z_1)-F(y_1)]}\right\}dz_1\biggr) dx \\
&=& \int_{y_1}^{z_2} \frac{f(x)F(y_1)}{F(x)(F(y_1)-F(x))}\biggr(-\log \frac{1-F(z_2)}{1-F(x)} + \log \left(1-[F(z_2)-F(x)]\right)\biggr) dx \\
&-& \int_{y_1}^{z_2} \frac{f(x)}{F(y_1)-F(x)}\biggr(-\log \frac{1-F(z_2)}{1-F(x)} + \log \frac{1-[F(z_2)-F(y_1)]}{1-[F(x)-F(y_1)]}\biggr) dx \\
&=& -\int_{y_1}^{z_2} \frac{f(x)F(y_1)}{F(x)(F(y_1)-F(x))} \log \frac{1-F(z_2)}{1-F(x)}dx \\
&+& \int_{y_1}^{z_2} \frac{f(x)F(y_1)}{F(x)(F(y_1)-F(x))} \log\left(1-[F(z_2)-F(x)]\right)dx \\
&+& \int_{y_1}^{z_2} \frac{f(x)}{F(y_1)-F(x)} \log \frac{1-F(z_2)}{1-F(x)}dx \\
&-& \int_{y_1}^{z_2} \frac{f(x)}{F(y_1)-F(x)} \log \frac{1-[F(z_2)-F(y_1)]}{1-[F(x)-F(y_1)]} dx.
\end{eqnarray*}

\begin{eqnarray*}
K_3 &=& \int \int \int_{D_x} \frac{f(x)f(y_2)f(z_1)}{(1-[F(z_1)-F(x)])(1-[F(z_1)-F(y_1)])F(y_1)} dx dy_2 dz_1 \\
&=&\frac{1}{F(y_1)} \int_{-\infty}^{y_1} f(y_2) dy_2 \int_{y_1}^{z_2} f(x)\\
&&\biggr(\int_{x}^{z_2} \frac{f(z_1)}{(1-[F(z_1)-F(x)])(1-[F(z_1)-F(y_1)])} dz_1 \biggr) dx  \\
&=& \int_{y_1}^{z_2} \frac{f(x)}{F(y_1)-F(x)} \biggr(\int_{x}^{z_2} \left\{\frac{f(z_1)}{1-[F(z_1)-F(x)]} - \frac{f(z_1)}{1-[F(z_1)-F(y_1)]}\right\} dz_1 \biggr) dx \\
&=& -\int_{y_1}^{z_2} \frac{f(x)}{F(y_1)-F(x)} \log\left(1-[F(z_2)-F(x)]\right)dx \\ 
&+& \int_{y_1}^{z_2} \frac{f(x)}{F(y_1)-F(x)} \log \frac{1-[F(z_2)-F(y_1)]}{1-[F(x)-F(y_1)]} dx. 
\end{eqnarray*}

\begin{eqnarray*}
K_4 &=& \int \int \int_{D_x} \frac{f(x)f(y_2)f(z_1)}{(1-[F(x)-F(y_1)])(1-[F(z_1)-F(y_1)])(1-F(z_1))} dx dy_2 dz_1 \\
&=& \int_{-\infty}^{y_1} f(y_2) dy_2 \int_{y_1}^{z_2} \frac{f(x)}{1-[F(x)-F(y_1)]}\\
&&\biggr(\int_{x}^{z_2} \frac{f(z_1)}{(1-[F(z_1)-F(y_1)])(1-F(z_1))} dz_1 \biggr) dx  \\
&=& \int_{y_1}^{z_2} \frac{f(x)}{1-[F(x)-F(y_1)]}\biggr(\int_{x}^{z_2}\left\{ \frac{f(z_1)}{1-F(z_1)} dz_1 - \frac{f(z_1)}{1-[F(z_1)-F(y_1)]}\right\} dz_1\biggr) dx  \\
&=& \int_{y_1}^{z_2} \frac{f(x)}{1-[F(x)-F(y_1)]}\biggr(-\log \frac{1-F(z_2)}{1-F(x)} + \log \frac{1-[F(z_2)-F(y_1)]}{1-[F(x)-F(y_1)]}\biggr) dx  \\
&=& - \int_{y_1}^{z_2} \frac{f(x)}{1-[F(x)-F(y_1)]} \log \frac{1-F(z_2)}{1-F(x)} dx \\
&+& \int_{y_1}^{z_2} \frac{f(x)}{1-[F(x)-F(y_1)]} \log \frac{1-[F(z_2)-F(y_1)]}{1-[F(x)-F(y_1)]} dx.  
\end{eqnarray*}

\begin{eqnarray*}
K_5 &=& \int \int \int_{D_x} \frac{f(x)f(y_2)f(z_1)}{(1-[F(x)-F(y_1)])(1-[F(z_1)-F(y_1)])F(y_1} dx dy_2 dz_1 \\
&=& \frac{1}{F(y_1)} \int_{-\infty}^{y_1} f(y_2) dy_2 \int_{y_1}^{z_2} \frac{f(x)}{1-[F(x)-F(y_1)]}\biggr(\int_{x}^{z_2} \frac{f(z_1)}{1-[F(z_1)-F(y_1)]} dz_1 \biggr) dx  \\
&=& - \int_{y_1}^{z_2} \frac{f(x)}{1-[F(x)-F(y_1)]} \log \frac{1-[F(z_2)-F(y_1)]}{1-[F(x)-F(y_1)]} dx.  
\end{eqnarray*}

\begin{eqnarray*}
K_6 &=& \int \int \int_{D_x} \frac{f(x)f(y_2)f(z_1)}{(1-[F(x)-F(y_1)])(1-F(x))(1-F(z_1)])} dx dy_2 dz_1 \\
&=& \int_{-\infty}^{y_1} f(y_2) dy_2 \int_{y_1}^{z_2} \frac{f(x)}{(1-F(x))(1-[F(x)-F(y_1))]}\biggr(\int_{x}^{z_2} \frac{f(z_1)}{1-F(z_1)} dz_1 \biggr) dx  \\
&=& - \int_{y_1}^{z_2} \frac{f(x)F(y_1)}{(1-F(x))(1-[F(x)-F(y_1))]} \log \frac{1-F(z_2)}{1-F(x)} dx.  
\end{eqnarray*}

\Bin Hence by regrouping all terms, we show that

\begin{eqnarray*}
&& K_1 + K_2 + K_3 + K_4 + K_5 + K_6 \\
&& = -\int_{y_1}^{z_2} \frac{f(x)}{F(x)(1-F(x))}\log \frac{1-F(z_2)}{1-F(x)} dx
\end{eqnarray*}

\Bin and this finishes the proof. $\blacksquare$ \\

\begin{corollary} \label{low_upp_ac_03}

The \textit{pdf} of $\left(Y_{(3)},Y^{(2)}\right)$ is then given by

$$
f_{(Y_{(3)},Y^{(2)})}(y,z) = -  f(y)f(z)\biggr(\int_{y}^{z} \frac{f(x)}{F(x)(1-F(x))}\log \frac{F(y)}{F(x)} dx\biggr)1_{(y<z)}.
$$ 

\end{corollary}

\Bin \textbf{Proof of Corollary \ref{low_upp_ac_03}}. We have 

\begin{eqnarray*}
f_{(Y_{(3)},Y^{(2)})}(y_2,z_1) &=& \int \int \int f_{Z_3}(x,y_1,y_2,z_1,z_2) dx dy_1 dz_2 \\
&=:& 1_{(y_2<z_1)} f(y_2)f(z_1)\left(M_1 + M_2 + M_3 + M_4 + M_5 + M_6\right).
\end{eqnarray*}

\Bin Let us denote by, for $(y_2,z_1)$ fixed,  $D_{(y_2,z_1)} = (y_2<y_1<x<z_1<z_2) =: D_x$ as the section of the domain of $Z_3$ at $(y_2,z_1)$ . So

\begin{eqnarray*}
M_1 &=& \int \int \int_{D_x} \frac{f(x)f(z_2)f(y_1)}{(1-[F(z_1)-F(x)])F(x)F(y_1)} dx dz_2 dy_1 \\
&=&\frac{1}{F(y_1)} \int_{z_1}^{+\infty} f(z_2) dz_2 \int_{y_2}^{z_1} \frac{f(x)}{F(x)(1-[F(z_1)-F(x)])}\\
&&\biggr(\int_{y_2}^{x} \frac{f(y_1)}{F(y_1)} dy_1 \biggr) dx  \\
&=& (1-F(z_1))\int_{y_2}^{z_1} \frac{f(x)}{F(x)(1-[F(z_1)-F(x)])} \log \frac{F(x)}{F(y_2)} dx.  
\end{eqnarray*}

\begin{eqnarray*}
M_2 &=& \int \int \int_{D_x} \frac{f(x)f(z_2)f(y_1)}{(1-[F(z_1)-F(x)])(1-[F(z_1)-F(y_1)])(1-F(z_1))} dx dz_2 dy_1 \\
&=& \frac{1}{1-F(z_1)} \int_{z_1}^{+\infty} f(z_2) dz_2 \int_{y_2}^{z_1} \frac{f(x)}{1-[F(z_1)-F(x)]}\\
&&\biggr(\int_{y_2}^{x} \frac{f(y_1)}{1-[F(z_1)-F(y_1)]} dy_1 \biggr) dx  \\
&=& \int_{y_2}^{z_1} \frac{f(x)}{F(x)(1-[F(z_1)-F(x)])} \log \frac{1-[F(z_1)-F(x)]}{1-[F(z_1)-F(y_2)]} dx.  
\end{eqnarray*}

\begin{eqnarray*}
M_3 &=& \int \int \int_{D_x} \frac{f(x)f(z_2)f(y_1)}{(1-[F(z_1)-F(x)])(1-[F(z_1)-F(y_1)])F(y_1)} dx dz_2 dy_1 \\
&=& \int_{z_1}^{+\infty} f(z_2) dz_2 \int_{y_2}^{z_1} \frac{f(x)}{1-[F(z_1)-F(x)]}\\
&&\biggr(\int_{y_2}^{x} \frac{f(y_1)}{(1-[F(z_1)-F(y_1))F(y_1)]} dy_1 \biggr) dx  \\
&=& \int_{y_2}^{z_1} \frac{f(x)}{1-[F(z_1)-F(x)]}\biggr(\int_{y_2}^{x} \left\{\frac{f(y_1)}{F(y_1)} - \frac{f(y_1)}{1-[F(z_1)-F(y_1)]} \right\}dy_1 \biggr) dx  \\
&=& \int_{y_2}^{z_1} \frac{f(x)}{1-[F(z_1)-F(x)]}\log \frac{F(x)}{F(y_2)} dx \\
&-& \int_{y_2}^{z_1} \frac{f(x)}{1-[F(z_1)-F(x)]}\log \frac{1-[F(z_1)-F(x)]}{1-[F(z_1)-F(y_2)]} dx.
\end{eqnarray*}

\begin{eqnarray*}
M_4 &=& \int \int \int_{D_x} \frac{f(x)f(z_2)f(y_1)}{(1-[F(x)-F(y_1)])(1-[F(z_1)-F(y_1)])(1-F(z_1))} dx dz_2 dy_1 \\
&=& \frac{1}{1-F(z_1)}\int_{z_1}^{+\infty} f(z_2) dz_2 \int_{y_2}^{z_1} f(x)\\
&&\biggr(\int_{y_2}^{x} \frac{f(y_1)}{(1-[F(x)-F(y_1)])(1-[F(z_1)-F(y_1)])} dy_1 \biggr) dx  \\
&=& \int_{y_2}^{z_1} \frac{f(x)}{F(z_1)-F(x)}\biggr(\int_{y_2}^{x} \left\{\frac{f(y_1)}{1-[F(z_1)-F(y_1)]} - \frac{f(y_1)}{1-[F(z_1)-F(y_1)]} \right\}dy_1 \biggr) dx  \\
&=& \int_{y_2}^{z_1} \frac{f(x)}{F(z_1)-F(x)} \log \frac{1-[F(z_1)-F(x)]}{1-[F(z_1)-F(y_2)]} dx \\
&-& \int_{y_2}^{z_1} \frac{f(x)}{F(z_1)-F(x)} \log \left(1-[F(x)-F(y_2)]\right) dx.
\end{eqnarray*}

\begin{eqnarray*}
M_5 &=& \int \int \int_{D_x} \frac{f(x)f(z_2)f(y_1)}{(1-[F(x)-F(y_1)])(1-[F(z_1)-F(y_1)])F(y_1)} dx dz_2 dy_1 \\
&=& \int_{z_1}^{+\infty} f(z_2) dz_2 \int_{y_2}^{z_1} f(x)\\
&&
\biggr(\int_{y_2}^{x} \frac{f(y_1)}{(1-[F(x)-F(y_1)])(1-[F(z_1)-F(y_1)])F(y_1)} dy_1 \biggr) dx.
\end{eqnarray*}

\Bin From now, we use the following relation

\begin{eqnarray*}
&&\frac{1}{(1-[F(x)-F(y_1)])(1-[F(z_1)-F(y_1)])F(y_1)} \\
&&=\frac{1}{(1-F(z_1))(F(z_1)-F(x))}\biggr(\frac{1}{F(y_1)} - \frac{1}{1-[F(z_1)-F(y_1)]}\biggr) \\
&&-\frac{1}{(1-F(x))(F(z_1)-F(x))}\biggr(\frac{1}{F(y_1)} - \frac{1}{1-[F(x)-F(y_1)]}\biggr).
\end{eqnarray*}
 
\Bin So we have 

\begin{eqnarray*}
M_5 &=& (1-F(z_1)) \int_{y_2}^{z_1} \frac{f(x)}{(1-F(z_1))(F(z_1)-F(x))}\\
&&\biggr(\int_{y_2}^{x} \left\{\frac{f(y_1)}{F(y_1)} - \frac{f(y_1)}{1-[F(z_1)-F(y_1)]} \right\}dy_1 \biggr) dx \\
&-& (1-F(z_1)) \int_{y_2}^{z_1} \frac{f(x)}{(1-F(z_1))(F(z_1)-F(x))}\\
&&\biggr(\int_{y_2}^{x} \left\{\frac{f(y_1)}{F(y_1)} - \frac{f(y_1)}{1-[F(x)-F(y_1)]} \right\}dy_1 \biggr) dx \\
&=& \int_{y_2}^{z_1} \frac{f(x)}{F(z_1)-F(x)}\biggr(\log \frac{F(x)}{F(y_2)} - \log \frac{1-[F(z_1)-F(x)]}{1-[F(z_1)-F(y_2)]} \biggr) dx  \\
&-& \int_{y_2}^{z_1} \frac{f(x)(1-F(z_1))}{(F(z_1)-F(x))(1-F(z_1))}\biggr(\log \frac{F(x)}{F(y_2)} + \log \left(1-[F(x)-F(y_2)]\right) \biggr) dx \\
&=& \int_{y_2}^{z_1} \frac{f(x)}{F(z_1)-F(x)} \log \frac{F(x)}{F(y_2)} dx \\
&-& \int_{y_2}^{z_1} \frac{f(x)}{F(z_1)-F(x)} \log \frac{1-[F(z_1)-F(x)]}{1-[F(z_1)-F(y_2)]} dx \\
&-& \int_{y_2}^{z_1} \frac{f(x)(1-F(z_1))}{(F(z_1)-F(x))(1-F(z_1))} \log \frac{F(x)}{F(y_2)} dx \\
&-& \int_{y_2}^{z_1} \frac{f(x)(1-F(z_1))}{(F(z_1)-F(x))(1-F(z_1))} \log \left(1-[F(x)-F(y_2)]\right) dx.
\end{eqnarray*}

\begin{eqnarray*}
M_6 &=& \int \int \int_{D_x} \frac{f(x)f(z_2)f(y_1)}{(1-[F(x)-F(y_1)])(1-F(x))(1-F(z_1))} dx dz_2 dy_1 \\
&=& (1-F(z_1))\int_{z_1}^{+\infty} f(z_2) dz_2 \int_{y_2}^{z_1} \frac{f(x)}{1-F(x)}\biggr(\int_{y_2}^{x} \frac{f(y_1)}{1-[F(x)-F(y_1)]} dy_1 \biggr) dx \\
&-& \int_{y_2}^{z_1} \frac{f(x)}{1-F(x)} \log \left(1-[F(x)-F(y_2)]\right) dx.
\end{eqnarray*}

\Bin Hence by regrouping all terms, we show that

\begin{eqnarray*}
&& M_1 + M_2 + M_3 + M_4 + M_5 + M_6 \\
&& = -\int_{y_2}^{z_1} \frac{f(x)}{F(x)(1-F(x))}\log \frac{F(y_2)}{F(x)} dx
\end{eqnarray*}

\Bin and this finishes the proof. $\blacksquare$ \\

\begin{corollary} \label{low_upp_ac_04}

\Bin The \textit{pdf} of $\left(Y_{(3)},Y^{(3)}\right)$ is then given by

$$
f_{(Y_{(3)},Y^{(3)})}(y,z) =  f(y)f(z)\biggr(\int_{y}^{z} \frac{f(x)}{F(x)(1-F(x))}\log \frac{F(y)}{F(x)}\log \frac{1-F(z)}{1-F(x)} dx\biggr)1_{(y<z)}.
$$
\end{corollary}

\Bin \textbf{Proof of Corollary \ref{low_upp_ac_04}}. We have 

\begin{eqnarray*}
f_{(Y_{(3)},Y^{(3)})}(y_2,z_2) &=& \int \int \int f_{Z_3}(x,y_1,y_2,z_1,z_2) dx dy_1 dz_1 \\
&=:& 1_{(y_2<z_2)} f(y_2)f(z_2)\left(N_1 + N_2 + N_3 + N_4 + N_5 + N_6\right).
\end{eqnarray*}

\Bin Let us denote by, for $(y_2,z_2)$ fixed,  $D_{(y_2,z_2)} = (y_2<y_1<x<z_1<z_2) =: D_x$ as the section of the domain of $Z_3$ at $(y_2,z_2)$ . So

\begin{eqnarray*}
N_1 &=& \int \int \int_{D_x} \frac{f(x)f(y_1)f(z_1)}{(1-[F(z_1)-F(x)])F(x)F(y_1)} dx dy_1 dz_1 \\
&=& \int_{y_2}^{z_2} \frac{f(x)}{F(x)}\biggr(\int_{x}^{z_2} \frac{f(z_1)}{(1-[F(z_1)-F(x)])} dz_1 \int_{y_2}^{x} \frac{f(y_1)}{F(y_1)}dy_1 \biggr) dx  \\
&=& -\int_{y_2}^{z_2} \frac{f(x) \log\left(1-[F(z_2)-F(x)]\right) \log \frac{F(x)}{F(y_2)}}{F(x)} dx.  
\end{eqnarray*}

\begin{eqnarray*}
N_2 &=& \int \int \int_{D_x} \frac{f(x)f(y_1)f(z_1)}{(1-[F(z_1)-F(x)])(1-[F(z_1)-F(y_1)])(1-F(z_1))} dx dy_1 dz_1 \\
&=& \int_{y_2}^{z_2} f(x)\biggr(\int_{x}^{z_2} \frac{f(z_1)}{(1-[F(z_1)-F(x)])(1-F(z_1))} \biggr(\int_{y_2}^{x} \frac{f(y_1)}{(1-[F(z_1)-F(y_1)])}dy_1\biggr) dz_1 \biggr) dx  \\
&=& \int_{y_2}^{z_2} f(x)\biggr(\int_{x}^{z_2}\frac{f(z_1) \log \frac{1-[F(z_1)-F(x)]}{1-[F(z_1)-F(y_2)]}}{(1-[F(z_1)-F(x)])(1-F(z_1))} dz_1 \biggr) dx.  
\end{eqnarray*}

\begin{eqnarray*}
N_3 &=& \int \int \int_{D_x} \frac{f(x)f(y_1)f(z_1)}{(1-[F(z_1)-F(x)])(1-[F(z_1)-F(y_1)])F(y_1)} dx dy_1 dz_1 \\
&=& \int_{y_2}^{z_2} f(x)\biggr(\int_{x}^{z_2} \frac{f(z_1)}{(1-[F(z_1)-F(x)])}\\
&& \biggr(\int_{y_2}^{x} \frac{f(y_1)}{(1-[F(z_1)-F(y_1)]F(y_1))}dy_1\biggr) dz_1 \biggr) dx  \\
&=& \int_{y_2}^{z_2} f(x)\biggr(\int_{x}^{z_2}\frac{f(z_1)}{(1-[F(z_1)-F(x)])(1-F(z_1))}\\
&& \biggr(\int_{y_2}^{x}\left\{\frac{f(y_1)}{F(y_1)}-\frac{f(y_1)}{(1-[F(z_1)-F(y_1)])}\right\} dy_1\biggr) dz_1 \biggr) dx \\
&=& \int_{y_2}^{z_2} f(x)\biggr(\int_{x}^{z_2}\frac{f(z_1)}{(1-[F(z_1)-F(x)])(1-F(z_1))}\\
&& \biggr(\log \frac{F(x)}{F(y_2)} - \log \frac{1-[F(z_1)-F(x)]}{1-[F(z_1)-F(y_2)]}\biggr) dz_1 \biggr) dx \\ 
&=& \int_{y_2}^{z_2} \frac{f(x)\log \frac{F(x)}{F(y_2)}}{F(x)}\\
&& \biggr(\int_{x}^{z_2}\left\{\frac{f(z_1)}{1-F(z_1)}-\frac{f(z_1)}{(1-[F(z_1)-F(x)])}\right\} dz_1 \biggr) dx - N_2 \\
&=& \int_{y_2}^{z_2} \frac{f(x)\log \frac{F(x)}{F(y_2)}}{F(x)} \biggr(-\log \frac{1-F(z_2)}{1-F(x)} + \log \left(1-[F(z_2)-F(x)]\right) \biggr) dx - N_2 \\
&=& -\int_{y_2}^{z_2} \frac{f(x)\log \frac{F(x)}{F(y_2)}\log \frac{1-F(z_2)}{1-F(x)}}{F(x)} dx - N_1 - N_2.
\end{eqnarray*}

\newpage
\begin{eqnarray*}
N_4 &=& \int \int \int_{D_x} \frac{f(x)f(y_1)f(z_1)}{(1-[F(x)-F(y_1)])(1-[F(z_1)-F(y_1)])(1-F(z_1))} dx dy_1 dz_1 \\
&=& \int_{y_2}^{z_2} f(x)\biggr(\int_{x}^{z_2}\frac{f(z_1)}{(1-F(z_1))}\\
&& \biggr(\int_{y_2}^{x} \frac{f(y_1)}{(1-[F(x)-F(y_1)])(1-[F(z_1)-F(y_1)])}dy_1\biggr) dz_1 \biggr) dx  \\
&=& \int_{y_2}^{z_2} f(x)\biggr(\int_{x}^{z_2}\frac{f(z_1)}{(1-F(z_1))(F(x)-F(z_1))} \biggr(\int_{y_2}^{x} \\
&&\biggr(\frac{f(y_1)}{(1-[F(x)-F(y_1)])} \\
&-& \frac{f(y_1)}{(1-[F(z_1)-F(y_1)])}\biggr)dy_1\biggr) dz_1 \biggr) dx  \\
&=& \int_{y_2}^{z_2} f(x)\biggr(\int_{x}^{z_2}\frac{f(z_1)}{(1-F(z_1))(F(x)-F(z_1))} \\
&&\biggr(-\log(1-[F(x)-F(y_2)]) \\
&-& \log \left(\frac{1-[F(z_1)-F(x)]}{1-[F(z_1)-F(y_2)]}\right)\biggr)dz_1 \biggr) dx \\
&=& -\int_{y_2}^{z_2} f(x)\log(1-[F(x)-F(y_2)])\biggr(\int_{x}^{z_2}\frac{f(z_1)}{(1-F(z_1))(F(x)-F(z_1))} dz_1 \biggr) dx \\
&-& \int_{y_2}^{z_2} f(x)\biggr(\int_{x}^{z_2}\frac{f(z_1)\log \left(\frac{1-[F(z_1)-F(x)]}{1-[F(z_1)-F(y_2)]}\right)}{(1-F(z_1))(F(x)-F(z_1))} dz_1 \biggr) dx.
\end{eqnarray*}

\Bin From now, we use the simple relation

\begin{eqnarray} \label{eq1}
&&\frac{1}{(1-F(z_1))(F(x)-F(z_1))} \\
&&=\frac{1}{(1-F(x))}\biggr(\frac{1}{F(x)-F(z_1)} - \frac{1}{1-F(z_1)}\biggr) 
\end{eqnarray}

\Bin to arrive at

\begin{eqnarray*}
N_4 &=& -\int_{y_2}^{z_2} \frac{f(x)\log(1-[F(x)-F(y_2)])}{1-F(x)}\biggr(\int_{x}^{z_2}\frac{f(z_1)}{F(x)-F(z_1)} dz_1 \biggr) dx \\
&+& \int_{y_2}^{z_2} \frac{f(x)\log(1-[F(x)-F(y_2)])}{1-F(x)}\biggr(\int_{x}^{z_2}\frac{f(z_1)}{1-F(z_1)} dz_1 \biggr) dx \\
&-& \int_{y_2}^{z_2} f(x)\biggr(\int_{x}^{z_2}\frac{f(z_1)\log \left(\frac{1-[F(z_1)-F(x)]}{1-[F(z_1)-F(y_2)]}\right)}{(1-F(z_1))(F(x)-F(z_1))} dz_1 \biggr) dx \\
&=& \int_{y_2}^{z_2} \frac{f(x)\log(1-[F(x)-F(y_2)])\log(F(x)-F(z_2))}{1-F(x)} dx \\
&-& \int_{y_2}^{z_2} \frac{f(x)\log(1-[F(x)-F(y_2)])\log \frac{1-F(z_2)}{1-F(x)}}{1-F(x)} dx \\
&-& \int_{y_2}^{z_2} f(x)\biggr(\int_{x}^{z_2}\frac{f(z_1)\log \left(\frac{1-[F(z_1)-F(x)]}{1-[F(z_1)-F(y_2)]}\right)}{(1-F(z_1))(F(x)-F(z_1))} dz_1 \biggr) dx.
\end{eqnarray*}

\begin{eqnarray*}
N_5 &=& \int \int \int_{D_x} \frac{f(x)f(y_1)f(z_1)}{(1-[F(x)-F(y_1)])(1-[F(z_1)-F(y_1)])F(y_1))} dx dy_1 dz_1 \\
&=& \int_{y_2}^{z_2} f(x)\biggr(\int_{x}^{z_2} f(z_1) \\
&&\biggr(\int_{y_2}^{x} \frac{f(y_1)}{(1-[F(x)-F(y_1)])(1-[F(z_1)-F(y_1)])F(y_1)}dy_1\biggr) dz_1 \biggr) dx  \\
&=& \int_{y_2}^{z_2} f(x)\biggr(\int_{x}^{z_2} \frac{f(z_1)}{(1-F(x))(F(x)-F(z_1))} \\
&&\biggr(\log \frac{F(x)}{F(y_2)} + \log \left(1-[F(x)-F(y_2)]\right)\biggr) dz_1 \biggr) dx  \\
&-& \int_{y_2}^{z_2} f(x)\biggr(\int_{x}^{z_2} \frac{f(z_1)}{(1-F(z_1))(F(x)-F(z_1))} \\
&&\biggr(\log \frac{F(x)}{F(y_2)} - \log \left(\frac{1-[F(z_1)-F(x)]}{1-[F(z_1)-F(y_2)]}\right)\biggr) dz_1 \biggr) dx   
\end{eqnarray*}

\Bin where in the last equality we used the simple following relation 

\begin{eqnarray*}
&&\frac{1}{(1-[F(x)-F(y_1)])(1-[F(z_1)-F(y_1)])F(y_1)} \\
&&=\frac{1}{(1-F(x))(F(x)-F(z_1))}\biggr(\frac{1}{F(y_1)} - \frac{1}{1-[F(x)-F(y_1)]}\biggr) \\
&&-\frac{1}{(1-F(z_1))(F(x)-F(z_1))}\biggr(\frac{1}{F(y_1)} - \frac{1}{1-[F(z_1)-F(y_1)]}\biggr).
\end{eqnarray*}
 
\Bin So we have 

\begin{eqnarray*}
N_5 &=& \int_{y_2}^{z_2} f(x)\biggr(\int_{x}^{z_2} \frac{f(z_1)\log \frac{F(x)}{F(y_2)}}{(1-F(x))(F(x)-F(z_1))}\biggr) dz_1 \biggr) dx \\
&+& \int_{y_2}^{z_2} f(x)\biggr(\int_{x}^{z_2} \frac{f(z_1)\log \left(1-[F(x)-F(y_2)]\right)}{(1-F(x))(F(x)-F(z_1))}\biggr) dz_1 \biggr) dx \\
&-& \int_{y_2}^{z_2} f(x)\biggr(\int_{x}^{z_2} \frac{f(z_1)\log \frac{F(x)}{F(y_2)}}{(1-F(z_1))(F(x)-F(z_1))}\biggr) dz_1 \biggr) dx \ \ (L3)\\
&+& \int_{y_2}^{z_2} f(x)\biggr(\int_{x}^{z_2} \frac{f(z_1)\log \left(\frac{1-[F(z_1)-F(x)]}{1-[F(z_1)-F(y_2)]}\right)}{(1-F(z_1))(F(x)-F(z_1))}\biggr) dz_1 \biggr) dx \\
&=& -\int_{y_2}^{z_2} \frac{f(x)\log \frac{F(x)}{F(y_2)}\log\left(F(x)-F(z_2)\right)}{1-F(x)} dx\\
&-& \int_{y_2}^{z_2} \frac{f(x)\log \left(1-[F(x)-F(y_2)]\right)\log\left(F(x)-F(z_2)\right)}{1-F(x)} dx\\  
&+& \int_{y_2}^{z_2} \frac{f(x)\log \frac{F(x)}{F(y_2)}\log\left(F(x)-F(z_2)\right)}{1-F(x)} dx\\
&-& \int_{y_2}^{z_2} \frac{f(x)\log \frac{F(x)}{F(y_2)}\log \left(\frac{1-F(z_2)}{1-F(x)}\right)}{1-F(x)} dx\\
&+& \int_{y_2}^{z_2} f(x)\biggr(\int_{x}^{z_2} \frac{f(z_1)\log \left(\frac{1-[F(z_1)-F(x)]}{1-[F(z_1)-F(y_2)]}\right)}{(1-F(z_1))(F(x)-F(z_1))}\biggr) dz_1 \biggr) dx
\end{eqnarray*}

\Bin where we use in Line (L3) we used the relation (\ref{eq1}) above

\begin{eqnarray*}
N_6 &=& \int \int \int_{D_x} \frac{f(x)f(y_1)f(z_1)}{(1-[F(x)-F(y_1)])(1-F(x))(1-F(z_1))} dx dy_1 dz_1 \\
&=& \int_{y_2}^{z_2} \frac{f(x)}{1-F(x)}\biggr( \int_{x}^{z_2} \frac{f(z_1)}{1-F(z_1)}dz_1 \int_{y_2}^{x} \frac{f(y_1)}{1-[F(x)-F(y_1)]}dy_1\biggr) dx  \\
&=& \int_{y_2}^{z_2} \frac{f(x)\log \left(\frac{1-F(z_2)}{1-F(x)}\right)\log \left(1-[F(x)-F(y_2)]\right)}{1-F(x)} dx.
\end{eqnarray*}

\Bin Hence, we have

$$
N_1 + N_2 + N_3 = -\int_{y_2}^{z_2} \frac{f(x)\log \frac{F(x)}{F(y_2)}\log \frac{1-F(z_2)}{1-F(x)}}{F(x)} dx
$$

\Bin and

$$
N_4 + N_5 + N_6 = -\int_{y_2}^{z_2} \frac{f(x)\log \frac{F(x)}{F(y_2)}\log \frac{1-F(z_2)}{1-F(x)}}{1-F(x)} dx
$$

\Bin and the proof is over by regrouping all terms. $\blacksquare$\\

\Ni Now, let us focus on discrete random variables.\\

\section{Discrete records} \label{sec_dis_records}

\noindent We will follow the steps we used in the above chapter by beginning by the case $n=2$. Here $f$ will denote the mass probability functions of the
random variables.\\
 
\Ni In our proofs below, we use the same graphical representations in Fig. \ref{fig1}, \ref{fig1b} and \ref{fig2}, in which the belongings
$(Y \in u^{\pm})$ are replaced by $(Y=u)$ with $Y$  playing the roles of the observation $Y_j$'s or the records values $Y^{(n)}$ or $Y_{(n)}$, and $u$ the roles
of the $y_n$'s or $z_n$'s.\\

\subsection{Probability law of the simultaneous joint lower-upper up to two records} \label{sec_01_ulrecords_chap3}  $ $\\

\Bin Below, we state the \textit{pdf}, give the proof and derive known results as means of verification of the results and get the \textit{pdf} of $(Y_{(p)},Y^{(q)})$ for $p\geq 2, \ q\geq 2$, when these records are defined. \\

\subsubsection{Finding the pdf of $Z_2$}

\begin{proposition} \label{prop_3} Let $Y_1$, $Y_2$, $\cdots$ be a sequence of independent and identically distributed discrete real random variables defined on the same probability space $\left(\Omega, \mathcal{A}, \mathbb{P}\right)$ with common discrete \textit{pdf} $f$ given on the strictly support 

$$
\nu_{Y} = \left\{y_j, \ j\in J\right\}, \ \ J\subset \mathbb{N},
$$

\Bin by 

$$
f(y) = \mathbb{P}(Y=y) \ \ \forall y\in \nu_{Y} \ \ and \ \ f(y)=0 \ \forall y\notin \nu_{Y}.
$$

\Bin The \textit{cdf} $F$ is given by

$$
F(y) = \mathbb{P}(Y\leq y) = \sum_{j\in J,\ y_j\leq y} f(y_j). 
$$

\Bin Let us denoted by

$$
F^{\ast}(y) = \mathbb{P}(Y<y) = \sum_{j\in J,\ y_j<y} f(y_j). 
$$

\Bin Suppose that $u(2)$ and $\ell(2)$ are finite. Then $Z_2$ has the following \textit{pdf}:

$$
f_{Z_2}(y_1,y_2,y_3)=\frac{f(y_1)f(y_2)f(y_3)}{F^{\ast}(y_1)(1-F(y_1))} 1_{\nu_{Y}(y_2<y_1<y_3)}.
$$
\end{proposition}

\Bin \textbf{Proof of Proposition \ref{prop_3}}. Let $y=(y_1,y_2,y_3)\in \nu_{Y}^3$ and let us put

$$
A(y):=\left(Y_1=y_1, Y_{(2)}=y_2,  Y^{(2)}=y_3\right).
$$

\Bin By definition, the \textit{pdf} of $Z_3$ is \\

$$
f_{Z_2}(y_1,y_2,y_3)=\mathbb{P}(A(y)), \ (y_1,y_2,y_3)\in \nu_{Y}^3.
$$

\Bin When dealing both with lower or upper records, and since the records are strong and an observation can be repeated with discrete random variables:\\

\Ni (a) the event $(u(2)=\ell(2))$ is negligible;\\

\Ni (b) $Y_2$ may be equal to $Y_1$ and so would not be the second record.\\

\Ni Now, we have 

$$
A(y)=\biggr(A(y) \cap (u(2)<\ell(2))\biggr) + \biggr(A(y) \cap (\ell(2)<u(2))\biggr). 
$$

\Bin Let $N_1$ be the inter-record time between the first and the second record time and let $N_2$ be the inter-record time between the two second record times. On $H_1:=(A(y) \cap (u(2)<\ell(2))$, as illustrated in Fig. \ref{fig1}, we have that $y_2<y_1<y_3$ and for

$$
G_0:=\biggr(Y_1=y_1, \ Y_2=y_3, \ Y_3=y_2\biggr),
$$

\Bin we have

$$
(\omega \in H_1)\cap (N_1=0,N_2=0) \Rightarrow (\omega \in G_{0}),
$$

\Bin and for $h\geq 1$, \ for $k\geq 1$,

\begin{eqnarray*}
(\omega \in H_1)\cap (N_1=h,N_2=k) &\Rightarrow& \biggr(\omega \in (Y_1=y_1) \bigcap (Y_{1+j}=y_1,\ 1\leq j\leq h)\bigcap(Y_{2+h}=y_3) \\
&&\bigcap (Y_{2+h+j}\geq y_1, \ 1\leq j\leq k) \bigcap ( Y_{3+h+k}=y_2)\biggr).
\end{eqnarray*}

\Bin Conversely, for $y_2<y_1<y_3$, for  

\begin{eqnarray*}
G_{(h,k)} &:=& \biggr(\omega \in (Y_1=y_1) \bigcap (Y_{1+j}=y_1,\ 1\leq j\leq h)\bigcap(Y_{2+h}=y_3)\bigcap (Y_{2+h+j}\geq y_1, \ 1\leq j\leq k)\\
&&\ \ \ \bigcap \ \ ( Y_{3+h+k}=y_2)\biggr), \ h\geq 1,\ k\geq 1,
\end{eqnarray*}

\Bin we have for any $h\geq 0$, \ $k\geq 0$,  $\omega \in G_{(h,k)}$ implies that $\omega \in A(y)\cap (N_1=h,N_2=k) \cap (u(2)<\ell(2))$ and hence

$$
G_{(h,k)}=H_1\cap (N_1=h,N_2=k), \ h\geq 0,\ k\geq 0.
$$

\Bin But

$$
\mathbb{P}(G_{(h,k)})= f(y_1)f(y_1)^h f(y_3) \mathbb{P}(Y\geq y_1)^k f(y_2) , \ h\geq 0, \ k\geq 0.
$$

\Bin Finally we have

$$
\mathbb{P}(H_1)=\sum_{h\geq 0} \sum_{k\geq 0} \mathbb{P}(G_{(h,k)}) = \frac{f(y_1)f(y_2)f(y_3)}{(1-f(y_1))F^{\ast}(y_1)} 1_{\nu_Y\cap (y_2<y_1<y_3)}. 
$$

\Bin We treat $(\omega \in H_2)\cap (N_1=h,N_2=k)$ in the same manner (see  Fig. \ref{fig1b}) where $H_2=(A(y) \cap (u(2)>\ell(2))$ and the $G_{(h,k)}$'s are replaced by, for $y_2<y_1<y_3$,

$$
G_0:=\biggr(Y_1=y_1, \ Y_2=y_2, \ Y_3=y_3\biggr),
$$

\Bin and

\begin{eqnarray*}
G_{(h,k)} &=& \biggr(\omega \in (Y_1=y_1) \bigcap (Y_{1+j}=y_1,\ 1\leq j\leq h)\bigcap(Y_{2+h}=y_2)\bigcap (Y_{2+h+j}\leq y_1, \ 1\leq j\leq k)\\
&&\ \ \ \bigcap \ \ ( Y_{3+h+k}=y_3)\biggr), \ h\geq 1,\ k\geq 1,
\end{eqnarray*}

\Bin and 

$$
\mathbb{P}(G_{(h,k)})= f(y_1)f(y_1)^h f(y_2) \mathbb{P}(Y\leq y_1)^k f(y_3) , \ h\geq 0, \ k\geq 0.
$$

\Bin We conclude that

$$
\mathbb{P}(H_2)=\sum_{h\geq 0} \sum_{k\geq 0} \mathbb{P}(G_{(h,k)}) = \frac{f(y_1)f(y_2)f(y_3)}{(1-f(y_1))(1-F(y_1))} 1_{\nu_Y\cap (y_2<y_1<y_3)}. 
$$

\Bin Hence 

\begin{eqnarray*}
\mathbb{P}(A(y)) &=& \mathbb{P}(H_1) + \mathbb{P}(H_2) \\
&=& \frac{f(y_1)f(y_2)f(y_3)}{(1-f(y_1))} \biggr(\frac{1}{F^{\ast}(y_1)} + \frac{1}{1-F(y_1)}\biggr) 1_{\nu_Y\cap (y_2<y_1<y_3)}\\
&=& \frac{f(y_1)f(y_2)f(y_3)}{F^{\ast}(y_1)(1-F(y_1))} 1_{\nu_Y\cap (y_2<y_1<y_3)}.
\end{eqnarray*}

\Bin The proof is over. $\blacksquare$\\

\subsubsection{Derivation of known results and of the second record values (lower and upper)} $ $ \\

\noindent \textbf{(i) Let us rediscover the law of $\left(Y_1,Y^{(2)}\right)$ whose \textit{pdf} is}:

\begin{eqnarray*}
f_{(Y_1,Y^{(2)})}(y_1,y_3) &=& \sum_{y_2} f_{Z_2}(y_1,y_2,y_3) \\
&=& 1_{(y_1<y_3)} \frac{f(y_1)f(y_3)}{F^{\ast}(y_1)(1-F(y_1))}  \sum_{y_2 \in \nu_{Y},\ y_2<y_1} f(y_2) \\
&=& \frac{f(y_1)}{1-F(y_1)} f(y_3)1_{(y_1<y_3)}.
\end{eqnarray*}

\Bin So, if $r(x)=f(x)/(1-F(x))$, $x\in ]\textit{\text{lep}}(F), \ \textit{\text{uep}}(F)[$, then

$$
f_{(Y_1,Y^{(2)})}(x,y) = r(x) f(y)1_{(x<y)}
$$

\Bin and this confirms a known result. \\

\noindent \textbf{(ii) We can derive the law of $\left(Y_{(2)},Y^{(2)}\right)$ as given below}:

\begin{corollary} \label{low_upp_d_01}

\Ni The \textit{pdf} of $\left(Y_{(2)},Y^{(2)}\right)$ is then 

$$
f_{(Y_{(2)},Y^{(2)})}(y,z) = f(y)f(z) \sum_{y<x<z} \frac{f(x)}{{F^{\ast}(x)(1-F(x))}} 1_{(y<z)}. 
$$

\end{corollary}

\Bin \textbf{Proof of Corollary \ref{low_upp_d_01}}. We have

\begin{eqnarray*}
f_{(Y_{(2)},Y^{(2)})}(y_2,y_3) &=& \sum_{y_1} f_{Z_2}(y_1,y_2,y_3) \\
&=& 1_{(y_2<y_3)}f(y_2)f(y_3) \sum_{y_2<y_1<y_3} \frac{f(y_1)}{{F^{\ast}(y_1)(1-F(y_1))}}
\end{eqnarray*}

\Bin and the proof is over. $\blacksquare$ \\

\subsection{Probability law of the simultaneous joint lower-upper up to two records} \label{sec_02_ulrecords_chap3}  $ $\\

\subsubsection{Finding the \textit{pdf} of $Z_3$}

\begin{proposition} \label{prop_4}
Let $Y_1$, $Y_2$, $\cdots$ be a sequence of independent and identically distributed discrete real random variables defined on the same probability space $\left(\Omega, \mathcal{A}, \mathbb{P}\right)$ with common discrete \textit{pdf} $f$ given on the strictly support 

$$
\nu_{Y} = \left\{y_j, \ j\in J\right\}, \ \ J\subset \mathbb{N},
$$

\Bin by 

$$
f(y) = \mathbb{P}(Y=y) \ \ \forall y\in \nu_{Y} \ \ and \ \ f(y)=0 \ \forall y\notin \nu_{Y}.
$$

\Bin The \textit{cdf} $F$ is given by

$$
F(y) = \mathbb{P}(Y\leq y) = \sum_{j\in J,\ y_j\leq y} f(y_j). 
$$

\Bin Let us denoted by

$$
F^{\ast}(y) = \mathbb{P}(Y<y) = \sum_{j\in J,\ y_j<y} f(y_j). 
$$

\Bin Suppose that $u(3)$ and $\ell(3)$ are finite. Then $Z_3$ has the following \textit{pdf}:

$$
f_{Z_3}(x,y_1,y_2,z_1,z_2)=f(x)f(y_1)f(y_2)f(z_1)f(z_2) L(x,y_1,z_1) 1_{\nu_{Y}\cap (y_2<y_1<x<z_1<z_2)},
$$

\Bin for

\begin{eqnarray*}
&&L(x,y_1,z_1)\\
&&=\frac{1}{(1-f(x))(1-[F(z_1)-F^{\ast}(x)])F^{\ast}(x)F^{\ast}(y_1)}\\
&&+\frac{1}{(1-f(x))(1-[F(z_1)-F^{\ast}(x)])(1-[F(z_1)-F^{\ast}(y_1)])(1-F(z_1))}\\
&&+\frac{1}{(1-f(x))(1-[F(z_1)-F^{\ast}(x)])(1-[F(z_1)-F^{\ast}(y_1)])F^{\ast}(y_1)}\\
&&+\frac{1}{(1-f(x))(1-[F(x)-F^{\ast}(y_1)])(1-[F(z_1)-F^{\ast}(y_1)])(1-F(z_1))}\\
&&+\frac{1}{(1-f(x))(1-[F(x)-F^{\ast}(y_1)])(1-[F(z_1)-F^{\ast}(y_1)])F^{\ast}(y_1)}\\
&&+\frac{1}{(1-f(x))(1-[F(x)-F^{\ast}(y_1)])(1-F(x))(1-F(z_1))}.
\end{eqnarray*}

\end{proposition}

\Bin \textbf{Proof of Proposition \ref{prop_4}}. We use similar methods to those in Section \ref{sec_01_ulrecords_chap3} just above but the situation a little more complex. We define

$$
A(y):=\left(Y_1=x, Y_{(2)}=y_1, Y_{(3)}=y_2, Y^{(2)}=z_1, Y^{(3)}=z_2\right),
$$

\Bin with $y=(x,y_1,y_2,z_1,z_2)\in \nu_{Y}^5$. We have to consider all the 24 orderings of
$(u(2),u(3),\ell(2),\ell(3))$. Fortunately $u(3)$ (resp. $\ell(3)$) cannot come before $u(2)$ (resp. $\ell(2)$). It will remains six orderings

\begin{eqnarray*}
O_1&=&(u(2)<u(3)<\ell(2)<\ell(3))\\
O_2&=&(u(2)<\ell(2)<\ell(3)<u(3))\\
O_3&=&(u(2)<\ell(2)<u(3)<\ell(3))\\
O_4&=&(\ell(2)<u(2)<\ell(3)<u(3))\\
O_5&=&(\ell(2)<u(2)<u(3)<\ell(3))\\
O_6&=&(\ell(2)<\ell(3)<u(2)<u(3)).
\end{eqnarray*} 

\Bin Let $N(h)$, $h \in \{1,2,3,4\}$ the four inter-record times. We have to decompose each $H_i=A(y)\cap O_i$, $i \in \{1,\cdots,6\}$ into

$$
H_i=\sum_{i_1\geq 0, i_2\geq 0, i_3\geq 0, i_4\geq 0} H_i \cap (N(1)=i_1,N(2)=i_2,N(3)=i_3,N(4)=i_4). 
$$

\Bin As above, we are going to describe $G_{i,i_1,i_2,i_3,i_4}=H_i \cap (N(1)=i_1,N(2)=i_2,N(3)=i_3,N(4)=i_4)$. We have six cases to deal with. So, we fully explain one cases and let the reader check the other five cases. From Fig. \ref{fig2}, we have the following facts :\\

\Ni (a) For $h \in [u(1),u(2)[$, $Y_h=x$ otherwise it would be the second record value (upper or lower);\\

\Ni (b) For $h \in [u(2),u(3)[$, $Y_h\geq x$ otherwise it would be the second lower record, and $Y_h\leq z_1$ otherwise it would be the third upper record;\\

\Ni (c) For $h \in [u(3),\ell(2)[$, $Y_h\geq x$ otherwise it would be the second lower record (but it is not bounded above).\\

\Ni (d) For $h \in [\ell(2),\ell(3)[$, $Y_h\geq y_1$ otherwise it would be the third lower record (but it is not bounded above).\\

\Ni Hence, we get 
$$
\mathbb{P}(G_{1,i_1,i_2,i_3,i_4}) = f(x)f(x)^{i_1}f(z_1)\mathbb{P}(x\leq Y\leq z_1)^{i_2}f(z_2)\mathbb{P}(Y\geq x)^{i_3}f(y_1)\mathbb{P}(Y\geq y_1)^{i_4}f(y_2),
$$

\Bin and hence

$$
\mathbb{P}(H_1) = \frac{f(x)f(y_1)f(y_2)f(z_1)f(z_2)}{(1-f(x))(1-[F(z_1)-F^{\ast}(x)])F^{\ast}(x)F^{\ast}(y_1)}.
$$

\Bin For the five other cases, we have

$$
\mathbb{P}(G_{2,i_1,i_2,i_3,i_4}) = f(x)f(y_1)f(y_2)f(z_1)f(z_2)f(x)^{i_1}\mathbb{P}(x\leq Y\leq z_1)^{i_2}\mathbb{P}(y_1\leq Y\leq z_1)^{i_3}\mathbb{P}(Y\leq z_1)^{i_4},
$$

$$
\mathbb{P}(G_{3,i_1,i_2,i_3,i_4}) = f(x)f(y_1)f(y_2)f(z_1)f(z_2)f(x)^{i_1}\mathbb{P}(x\leq Y\leq z_1)^{i_2}\mathbb{P}(y_1\leq Y\leq z_1)^{i_3}\mathbb{P}(Y\geq y_1)^{i_4},
$$

$$
\mathbb{P}(G_{4,i_1,i_2,i_3,i_4}) = f(x)f(y_1)f(y_2)f(z_1)f(z_2)f(x)^{i_1}\mathbb{P}(y_1\leq Y\leq x)^{i_2}\mathbb{P}(y_1\leq Y\leq z_1)^{i_3}\mathbb{P}(Y\leq z_1)^{i_4},
$$

$$
\mathbb{P}(G_{5,i_1,i_2,i_3,i_4}) = f(x)f(y_1)f(y_2)f(z_1)f(z_2)f(x)^{i_1}\mathbb{P}(y_1\leq Y\leq x)^{i_2}\mathbb{P}(y_1\leq Y\leq z_1)^{i_3}\mathbb{P}(Y\geq y_1)^{i_4},
$$

$$
\mathbb{P}(G_{6,i_1,i_2,i_3,i_4}) = f(x)f(y_1)f(y_2)f(z_1)f(z_2)f(x)^{i_1}\mathbb{P}(y_1\leq Y\leq x)^{i_2}\mathbb{P}(Y\leq x)^{i_3}\mathbb{P}(Y\leq z_1)^{i_4}.
$$

\Bin These facts when put together lead to

$$
f_{Z_3}(x,y_1,y_2,z_1,z_2)= f(x)f(y_1)f(y_2)f(z_1)f(z_2) L(x,y_1,z_1)1_{(y_2<y_1<x<z_1<z_2)}
$$

\Bin with

\begin{eqnarray*}
&&L(x,y_1,z_1)\\
&&=\frac{1}{(1-f(x))(1-[F(z_1)-F^{\ast}(x)])F^{\ast}(x)F^{\ast}(y_1)}\\
&&+\frac{1}{(1-f(x))(1-[F(z_1)-F^{\ast}(x)])(1-[F(z_1)-F^{\ast}(y_1)])(1-F(z_1))}\\
&&+\frac{1}{(1-f(x))(1-[F(z_1)-F^{\ast}(x)])(1-[F(z_1)-F^{\ast}(y_1)])F^{\ast}(y_1)}\\
&&+\frac{1}{(1-f(x))(1-[F(x)-F^{\ast}(y_1)])(1-[F(z_1)-F^{\ast}(y_1)])(1-F(z_1))}\\
&&+\frac{1}{(1-f(x))(1-[F(x)-F^{\ast}(y_1)])(1-[F(z_1)-F^{\ast}(y_1)])F^{\ast}(y_1)}\\
&&+\frac{1}{(1-f(x))(1-[F(x)-F^{\ast}(y_1)])(1-F(x))(1-F(z_1))}.
\end{eqnarray*}

\Bin

\subsubsection{Derivation of known results and of the third record values (lower and upper)} $ $\\

\noindent \textbf{(i) Checking of results by comparing with the known \textit{pdf} of $\left(Y_1, Y^{(2)}, Y^{(3)}\right)$}.\\

\Ni We have 

\begin{eqnarray*}
f_{(Y_1,Y^{(2)},Y^{(3)})}(x,z_1,z_2) &=& \sum_{y_1} \sum_{y_2} f_{Z_3}(x,y_1,y_2,z_1,z_2) \\
&=:& 1_{(x<z_1<z_2)} f(x)f(z_1)f(z_2)\left(S_1 + S_2 + S_3 + S_4 + S_5 + S_6\right).
\end{eqnarray*}

\Ni Let us denote, for $(x,z_1,z_2)$ fixed, by $D_x := D_{(x,z_1,z_2)} = (y_2<y_1<x)$ as the section of the domain of $Z_3$ at $(x,z_1,z_2)$. So

\begin{eqnarray*}
S_1 &=&  \sum_{y_1,y_2,\ (y_1,y_2)\in D_x} \frac{f(y_1)f(y_2)}{(1-f(x))(1-[F(z_1)-F^{\ast}(x)])F^{\ast}(x)F^{\ast}(y_1)} \\
&=& \frac{1}{(1-f(x))(1-[F(z_1)-F^{\ast}(x)])F^{\ast}(x)} \sum_{y_1<x} \frac{f(y_1)}{F^{\ast}(y_1)} \sum_{y_2<y_1} f(y_2)  \\
&=& \frac{1}{(1-f(x))(1-[F(z_1)-F^{\ast}(x)])}.
\end{eqnarray*}

\begin{eqnarray*}
S_2 &=& \sum_{y_1,y_2,\ (y_1,y_2)\in D_x} \frac{f(y_1)f(y_2)}{(1-f(x))(1-[F(z_1)-F^{\ast}(x)])(1-[F(z_1)-F^{\ast}(y_1)])(1-F(z_1))} \\
&=& \frac{1}{(1-f(x))(1-[F(z_1)-F^{\ast}(x)])(1-F(z_1))} \sum_{y_1<x} \frac{f(y_1)}{1-[F(z_1)-F^{\ast}(y_1)]} \sum_{y_2<y_1} f(y_2)  \\
&=& \frac{1}{(1-f(x))(1-[F(z_1)-F^{\ast}(x)])(1-F(z_1))} \sum_{y_1<x} \frac{f(y_1) F^{\ast}(y_1)}{1-[F(z_1)-F^{\ast}(y_1)]}.
\end{eqnarray*}

\Bin But remark that 

$$
\sum_{y_1<x} \frac{f(y_1) F^{\ast}(y_1)}{1-[F(z_1)-F^{\ast}(y_1)]} = F^{\ast}(x) - (1-F(z_1))\sum_{y_1<x} \frac{f(y_1)}{1-[F(z_1)-F^{\ast}(y_1)]}
$$

\Bin to have

\begin{eqnarray*}
S_2 &=& \frac{F^{\ast}(x)}{(1-f(x))(1-[F(z_1)-F^{\ast}(x)])(1-F(z_1))} \\
&-&  \frac{1}{(1-f(x))(1-[F(z_1)-F^{\ast}(x)])} \sum_{y_1<x} \frac{f(y_1)}{1-[F(z_1)-F^{\ast}(y_1)]}.
\end{eqnarray*}

\begin{eqnarray*}
S_3 &=& \sum_{y_1,y_2,\ (y_1,y_2)\in D_x} \frac{f(y_1)f(y_2)}{(1-f(x))(1-[F(z_1)-F^{\ast}(x)])(1-[F(z_1)-F^{\ast}(y_1)])F^{\ast}(y_1)} \\
&=& \frac{1}{(1-f(x))(1-[F(z_1)-F^{\ast}(x)])} \sum_{y_1<x} \frac{f(y_1)}{(1-[F(z_1)-F^{\ast}(y_1)])F^{\ast}(y_1)} \sum_{y_2<y_1} f(y_2)  \\
&=& \frac{1}{(1-f(x))(1-[F(z_1)-F^{\ast}(x)])} \sum_{y_1<x} \frac{f(y_1)}{1-[F(z_1)-F^{\ast}(y_1)]}.
\end{eqnarray*}

\begin{eqnarray*}
S_4 &=& \sum_{y_1,y_2,\ (y_1,y_2)\in D_x} \frac{f(y_1)f(y_2)}{(1-f(x))(1-[F(x)-F^{\ast}(y_1)])(1-[F(z_1)-F^{\ast}(y_1)])(1-F(z_1))} \\
&=& \frac{1}{(1-f(x))(1-F(z_1))} \sum_{y_1<x} \frac{f(y_1)}{(1-[F(x)-F^{\ast}(y_1)])(1-[F(z_1)-F^{\ast}(y_1)])} \sum_{y_2<y_1} f(y_2)  \\
&=& \frac{1}{(1-f(x))(1-F(z_1))} \sum_{y_1<x} \frac{f(y_1)F^{\ast}(y_1)}{(1-[F(x)-F^{\ast}(y_1)])(1-[F(z_1)-F^{\ast}(y_1)])}.
\end{eqnarray*}

\Bin From now, we use the following decomposition 

\begin{eqnarray*}
&&\frac{1}{(1-[F(x)-F^{\ast}(y_1)])(1-[F(z_1)-F^{\ast}(y_1)])} \\ 
&&=\frac{1}{F(z_1)-F(x)} \biggr( \frac{1}{1-[F(z_1)-F^{\ast}(y_1)]} - \frac{1}{1-[F(x)-F^{\ast}(y_1)}\biggr)
\end{eqnarray*}

\Bin to get

\begin{eqnarray*}
S_4 &=&  \frac{1}{(1-f(x))(1-F(z_1))(F(z_1)-F(x))} \sum_{y_1<x} \frac{f(y_1)F^{\ast}(y_1)}{1-[F(z_1)-F^{\ast}(y_1)]} \\
&-& \frac{1}{(1-f(x))(1-F(z_1))(F(z_1)-F(x))} \sum_{y_1<x} \frac{f(y_1)F^{\ast}(y_1)}{1-[F(x)-F^{\ast}(y_1)]}.
\end{eqnarray*}

\Bin But remark that 

$$
\sum_{y_1<x} \frac{f(y_1) F^{\ast}(y_1)}{1-[F(z_1)-F^{\ast}(y_1)]} = F^{\ast}(x) - (1-F(z_1))\sum_{y_1<x} \frac{f(y_1)}{1-[F(z_1)-F^{\ast}(y_1)]}
$$

\Bin and

$$
\sum_{y_1<x} \frac{f(y_1) F^{\ast}(y_1)}{1-[F(x)-F^{\ast}(y_1)]} = F^{\ast}(x) - (1-F(x))\sum_{y_1<x} \frac{f(y_1)}{1-[F(x)-F^{\ast}(y_1)]}
$$

\Bin and so

\begin{eqnarray*}
S_4 &=&  \frac{F^{\ast}(x)}{(1-f(x))(1-F(z_1))(F(z_1)-F(x))} \\ 
&-& \frac{1}{(1-f(x))(F(z_1)-F(x))} \sum_{y_1<x} \frac{f(y_1)F^{\ast}(y_1)}{1-[F(z_1)-F^{\ast}(y_1)]} \\
&-& \frac{F^{\ast}(x)}{(1-f(x))(1-F(z_1))(F(z_1)-F(x))} \\ 
&+& \frac{1-F(x)}{(1-f(x))(1-F(z_1))(F(z_1)-F(x))}\sum_{y_1<x} \frac{f(y_1)F^{\ast}(y_1)}{1-[F(x)-F^{\ast}(y_1)]}.
\end{eqnarray*}

\begin{eqnarray*}
S_5 &=& \sum_{y_1,y_2,\ (y_1,y_2)\in D_x} \frac{f(y_1)f(y_2)}{(1-f(x))(1-[F(x)-F^{\ast}(y_1)])(1-[F(z_1)-F^{\ast}(y_1)])F^{\ast}(y_1)} \\
&=& \frac{1}{1-f(x)} \sum_{y_1<x} \frac{f(y_1)}{(1-[F(x)-F^{\ast}(y_1)])(1-[F(z_1)-F^{\ast}(y_1)])F^{\ast}(y_1)} \sum_{y_2<y_1} f(y_2)  \\
&=& \frac{1}{1-f(x)} \sum_{y_1<x} \frac{f(y_1)}{(1-[F(x)-F^{\ast}(y_1)])(1-[F(z_1)-F^{\ast}(y_1)])}. 
\end{eqnarray*}

\Bin From now, we use the decomposition above to arrive at 

\begin{eqnarray*}
S_5 &=& \frac{1}{(1-f(x))(F(z_1)-F(x))} \sum_{y_1<x} \frac{f(y_1)}{1-[F(z_1)-F^{\ast}(y_1)]} \\
&-& \frac{1}{(1-f(x))(F(z_1)-F(x))} \sum_{y_1<x} \frac{f(y_1)}{1-[F(x)-F^{\ast}(y_1)]}
\end{eqnarray*}

\begin{eqnarray*}
S_6 &=& \sum_{y_1,y_2,\ (y_1,y_2)\in D_x} \frac{f(y_1)f(y_2)}{(1-f(x))(1-[F(x)-F^{\ast}(y_1)])(1-F(x))(1-F(z_1))} \\
&=& \frac{1}{(1-f(x))(1-F(x))(1-F(z_1))} \sum_{y_1<x} \frac{f(y_1)}{1-[F(x)-F^{\ast}(y_1)]} \sum_{y_2<y_1} f(y_2)  \\
&=& \frac{1}{(1-f(x))(1-F(x))(1-F(z_1))} \sum_{y_1<x} \frac{f(y_1)F^{\ast}(y_1)}{1-[F(x)-F^{\ast}(y_1)]}. 
\end{eqnarray*}

\Bin By using the remark above, we will have

\begin{eqnarray*}
S_6 &=& \frac{F^{\ast}(x)}{(1-f(x))(1-F(x))(1-F(z_1))} - \frac{1}{(1-f(x))(1-F(z_1))}\sum_{y_1<x} \frac{f(y_1)}{1-[F(x)-F^{\ast}(y_1)]}. 
\end{eqnarray*}

\Bin Hence by doing the simple computations, all terms together give

$$
S_1 + S_2 + S_3 + S_4 + S_5 + S_6 = \frac{1}{(1-F(x))(1-F(z_1))}
$$

\Bin and hence

\begin{eqnarray*}
f_{(Y_1,Y^{(2)},Y^{(3)})}(x,z_1,z_2) &=&  \frac{f(x)f(z_1)f(z_2)}{(1-F(x))(1-F(z_1))} 1_{(x<z_1<z_2)} \\
&=& r(x) r(z_1) f(z_2) 1_{(x<z_1<z_2)} 
\end{eqnarray*}

\Bin and this confirms a known result, with

$$
r(t)=f(t)/(1-F(t)), \ \  t\in ]\textit{\text{lep}}(F), \ \textit{\text{uep}}(F)[.
$$ 

\Bin  \textbf{(ii) Let us derive the laws of $\left(Y_{(p)},Y^{(q)}\right)$, $2\leq p,q \leq 3$}.\\

\Ni \textbf{Let us begin by checking the law of \textbf{$\left(Y_{(2)},Y^{(2)}\right)$}, which has already been found in Corollary \ref{low_upp_d_01}.} \\

\Ni We have

\begin{eqnarray*}
f_{(Y_{(2)},Y^{(2)})}(y_1,z_1) &=& \sum_{x} \sum_{y_2} \sum_{z_2} f_{Z_3}(x,y_1,y_2,z_1,z_2) \\
&=:& 1_{(y_1<z_1)} f(y_1)f(z_1)\left(A_1 + A_2 + A_3 + A_4 + A_5 + A_6\right).
\end{eqnarray*}

\Bin Let us denote by, for $(y_1,z_1)$ fixed,  $D_{(y_1,z_1)} = (y_2<y_1<x<z_1<z_2) =: D_x$ as the section of the domain of $Z_3$ at $(y_1,z_1)$ . So

\begin{eqnarray*}
A_1 &=&  \sum_{x,y_2,z_2,\ (x,y_2,z_2)\in D_x} \frac{f(x)f(y_2)f(z_2)}{(1-f(x))(1-[F(z_1)-F^{\ast}(x)])F^{\ast}(x)F^{\ast}(y_1)} \\
&=& \frac{1}{F^{\ast}(y_1)} \sum_{y_1<x<z_1} \frac{f(x)}{(1-f(x))(1-[F(z_1)-F^{\ast}(x)])F^{\ast}(x)} \sum_{y_2<y_1} f(y_2)\sum_{z_2>z_1} f(z_2)   \\
&=& \sum_{y_1<x<z_1} \frac{(1-F(z_1))f(x)}{(1-f(x))(1-[F(z_1)-F^{\ast}(x)])F^{\ast}(x)}.
\end{eqnarray*}

\Bin But remark that 

\begin{eqnarray*}
&&\frac{1-F(z_1)}{(1-f(x))(1-[F(z_1)-F^{\ast}(x)])F^{\ast}(x)} \\ 
&&=\frac{1}{F^{\ast}(x)} - \frac{1}{1-[F(z_1)-F^{\ast}(x)]}
\end{eqnarray*}

\Bin to have

\begin{eqnarray*}
A_1 &=& \sum_{y_1<x<z_1} \frac{f(x)}{(1-f(x))F^{\ast}(x)} \\
&-& \sum_{y_1<x<z_1} \frac{f(x)}{(1-f(x))(1-[F(z_1)-F^{\ast}(x)])}.
\end{eqnarray*}

\begin{eqnarray*}
&&A_2 \\
&&=  \sum_{x,y_2,z_2,\ (x,y_2,z_2)\in D_x} \frac{f(x)f(y_2)f(z_2)}{(1-f(x))(1-[F(z_1)-F^{\ast}(x)])(1-[F(z_1)-F^{\ast}(y_1)])(1-F(z_1))} \\
&&= \frac{1}{(1-[F(z_1)-F^{\ast}(y_1)])(1-F(z_1))} \sum_{y_1<x<z_1} \frac{f(x)}{(1-f(x))(1-[F(z_1)-F^{\ast}(x)])}\\
&& \sum_{y_2<y_1} f(y_2)\sum_{z_2>z_1} f(z_2)   \\
&&= \frac{F^{\ast}(y_1)}{1-[F(z_1)-F^{\ast}(y_1)]} \sum_{y_1<x<z_1} \frac{f(x)}{(1-f(x))(1-[F(z_1)-F^{\ast}(x)])}.
\end{eqnarray*}

\begin{eqnarray*}
A_3 &=&  \sum_{x,y_2,z_2,\ (x,y_2,z_2)\in D_x} \frac{f(x)f(y_2)f(z_2)}{(1-f(x))(1-[F(z_1)-F^{\ast}(x)])(1-[F(z_1)-F^{\ast}(y_1)])F^{\ast}(y_1)} \\
&=& \frac{1}{(1-[F(z_1)-F^{\ast}(y_1)])F^{\ast}(y_1)} \sum_{y_1<x<z_1} \frac{f(x)}{(1-f(x))(1-[F(z_1)-F^{\ast}(x)])}\\
&& \sum_{y_2<y_1} f(y_2)\sum_{z_2>z_1} f(z_2)   \\
&=& \frac{1-F(z_1)}{1-[F(z_1)-F^{\ast}(y_1)]} \sum_{y_1<x<z_1} \frac{f(x)}{(1-f(x))(1-[F(z_1)-F^{\ast}(x)])}.
\end{eqnarray*}

\begin{eqnarray*}
A_4 &=& \sum_{x,y_2,z_2,\ (x,y_2,z_2)\in D_x} \frac{f(x)f(y_2)f(z_2)}{(1-f(x))(1-[F(x)-F^{\ast}(y_1)])(1-[F(z_1)-F^{\ast}(y_1)])(1-F(z_1))} \\
&=& \frac{1}{(1-[F(z_1)-F^{\ast}(y_1)])(1-F(z_1))} \sum_{y_1<x<z_1} \frac{f(x)}{(1-f(x))(1-[F(x)-F^{\ast}(y_1)]) }\\
&\times& \sum_{y_2<y_1} f(y_2)\sum_{z_2>z_1} f(z_2)   \\
&=& \frac{F^{\ast}(y_1)}{1-[F(z_1)-F^{\ast}(y_1)]} \sum_{y_1<x<z_1} \frac{f(x)}{(1-f(x))(1-[F(z_1)-F^{\ast}(x)])}.
\end{eqnarray*}

\begin{eqnarray*}
A_5 &=& \sum_{x,y_2,z_2,\ (x,y_2,z_2)\in D_x} \frac{f(x)f(y_2)f(z_2)}{(1-f(x))(1-[F(x)-F^{\ast}(y_1)])(1-[F(z_1)-F^{\ast}(y_1)])F^{\ast}(y_1)} \\
&=& \frac{1}{(1-[F(z_1)-F^{\ast}(y_1)])F^{\ast}(y_1)} \sum_{y_1<x<z_1} \frac{f(x)}{(1-f(x))(1-[F(x)-F^{\ast}(y_1)])}\\
&& \sum_{y_2<y_1} f(y_2)\sum_{z_2>z_1} f(z_2)   \\
&=& \frac{1-F(z_1)}{1-[F(z_1)-F^{\ast}(y_1)]} \sum_{y_1<x<z_1} \frac{f(x)}{(1-f(x))(1-[F(z_1)-F^{\ast}(x)])}.
\end{eqnarray*}

\begin{eqnarray*}
A_6 &=& \sum_{x,y_2,z_2,\ (x,y_2,z_2)\in D_x} \frac{f(x)f(y_2)f(z_2)}{(1-f(x))(1-[F(x)-F^{\ast}(y_1)])(1-F(x))(1-F(z_1))} \\
&=& \frac{1}{1-F(z_1)} \sum_{y_1<x<z_1} \frac{f(x)}{(1-f(x))(1-[F(x)-F^{\ast}(y_1)])(1-F(x))}\\
&& \sum_{y_2<y_1} f(y_2)\sum_{z_2>z_1} f(z_2)   \\
&=& F^{\ast}(y_1) \sum_{y_1<x<z_1} \frac{f(x)}{(1-f(x))(1-[F(x)-F^{\ast}(y_1)])(1-F(x))}.
\end{eqnarray*}

\Bin Hence by doing direct computations, we will have

$$
A_1 + A_2 + A_3 = \sum_{y_1<x<z_1} \frac{f(x)}{(1-f(x))F^{\ast}(x)}
$$

\Bin and 

$$
A_4 + A_5 + A_6 = \sum_{y_1<x<z_1} \frac{f(x)}{(1-f(x))(1-F(x))}
$$

\Bin So by regrouping all terms together, we will have

$$
f_{(Y_{(2)},Y^{(2)})}(y_1,z_1) = f(y_1)f(z_1) \biggr(\sum_{y_1<x<z_1} \frac{f(x)}{F^{\ast}(x)(1-F(x))}\biggr)1_{(y_1<z_1)}
$$

\Bin and that confirms the result found in Corollary \ref{low_upp_d_01}. \\

\Ni \textbf{Now, we are going to derive the laws of $\left(Y_{(2)},Y^{(3)}\right)$, $\left(Y_{(3)},Y^{(2)}\right)$ and $\left(Y_{(3)},Y^{(3)}\right)$ in the next Corollaries.} \\

\begin{corollary} \label{low_upp_d_02}
The \textit{pdf} of $\left(Y_{(2)},Y^{(3)}\right)$ is then given by

$$
f_{(Y_{(2)},Y^{(3)})}(y_1,z_2) = 1_{(y_1<z_2)}f(y_1)f(z_2) \sum_{y_1<x<z_2} \frac{f(x)}{F^{\ast}(x)(1-F(x))} \sum_{x<z_1<z_2} \frac{f(z_1)}{1-F(z_1)}.
$$

\end{corollary}

\Bin \textbf{Proof of Corollary \ref{low_upp_d_02}}.
We have 

\begin{eqnarray*}
f_{(Y_{(2)},Y^{(3)})}(y_1,z_2) &=& \sum_{x} \sum_{y_2} \sum_{z_1} f_{Z_3}(x,y_1,y_2,z_1,z_2) \\
&=:& 1_{(y_1<z_2)} f(y_1)f(z_2)\left(B_1 + B_2 + B_3 + B_4 + B_5 + B_6\right).
\end{eqnarray*}

\Bin Let us denote, for $(y_1,z_2)$ fixed, by  $D_{(y_1,z_2)} = (y_2<y_1<x<z_1<z_2) =: D_x$ as the section of the domain of $Z_3$ at $(y_1,z_2)$. So

\begin{eqnarray*}
B_1 &=&  \sum_{x,y_2,z_1,\ (x,y_2,z_1)\in D_x} \frac{f(x)f(y_2)f(z_1)}{(1-f(x))(1-[F(z_1)-F^{\ast}(x)])F^{\ast}(x)F^{\ast}(y_1)} \\
&=& \frac{1}{F^{\ast}(y_1)} \sum_{y_1<x<z_2} \frac{f(x)}{(1-f(x))F^{\ast}(x)} \sum_{x<z_1<z_2} \frac{f(z_1)}{1-[F(z_1)-F^{\ast}(x)]} \sum_{y_2<y_1} f(y_2)   \\
&=& \sum_{y_1<x<z_2} \frac{f(x)}{(1-f(x))F^{\ast}(x)} \sum_{x<z_1<z_2} \frac{f(z_1)}{1-[F(z_1)-F^{\ast}(x)]}.
\end{eqnarray*}

\begin{eqnarray*}
B_2 &=&  \sum_{x,y_2,z_1,\ (x,y_2,z_1)\in D_x} \frac{f(x)f(y_2)f(z_1)}{(1-f(x))(1-[F(z_1)-F^{\ast}(x)])(1-[F(z_1)-F^{\ast}(y_1)])(1-F(z_1))} \\
&=& \sum_{y_1<x<z_2} \frac{f(x)}{1-f(x)} \sum_{x<z_1<z_2} \frac{f(z_1)}{(1-[F(z_1)-F^{\ast}(x)])(1-[F(z_1)-F^{\ast}(y_1)])(1-F(z_1))}\\
&& \sum_{y_2<y_1} f(y_2)   \\
&=& F^{\ast}(y_1) \sum_{y_1<x<z_2}\frac{f(x)}{1-f(x)} \sum_{x<z_1<z_2} \frac{f(z_1)}{(1-[F(z_1)-F^{\ast}(x)])(1-[F(z_1)-F^{\ast}(y_1)])(1-F(z_1))}.
\end{eqnarray*}

\Bin From now, we use the following decomposition

\begin{eqnarray*}
&&\frac{1}{(1-[F(z_1)-F^{\ast}(x)])(1-[F(z_1)-F^{\ast}(y_1)])(1-F(z_1))} \\
&&=\frac{1}{F^{\ast}(x)F^{\ast}(y_1)}\biggr(\frac{1}{1-F(z_1)} - \frac{1}{1-[F(z_1)-F^{\ast}(x)]}\biggr) \\
&&-\frac{1}{F^{\ast}(x)(F^{\ast}(x)-F^{\ast}(y_1))}\biggr(\frac{1}{1-[F(z_1)-F^{\ast}(y_1)])} - \frac{1}{1-[F(z_1)-F^{\ast}(x)]}\biggr)
\end{eqnarray*}

\Bin to get

\begin{eqnarray*}
B_2 &=& \sum_{y_1<x<z_2}\frac{f(x)}{F^{\ast}(x)(1-f(x))} \sum_{x<z_1<z_2} \frac{f(z_1)}{1-F(z_1)} \\
&-& \sum_{y_1<x<z_2}\frac{f(x)}{F^{\ast}(x)(1-f(x))} \sum_{x<z_1<z_2} \frac{f(z_1)}{1-[F(z_1)-F^{\ast}(y_1)]} \\
&-& F^{\ast}(y_1)\sum_{y_1<x<z_2}\frac{f(x)}{F^{\ast}(x)(F^{\ast}(x)-F^{\ast}(y_1))(1-f(x))} \sum_{x<z_1<z_2} \frac{f(z_1)}{1-[F(z_1)-F^{\ast}(y_1)]} \\
&+& F^{\ast}(y_1)\sum_{y_1<x<z_2}\frac{f(x)}{F^{\ast}(x)(F^{\ast}(x)-F^{\ast}(y_1))(1-f(x))} \sum_{x<z_1<z_2} \frac{f(z_1)}{1-[F(z_1)-F^{\ast}(x)]}. 
\end{eqnarray*}

\begin{eqnarray*}
B_3 &=&  \sum_{x,y_2,z_1,\ (x,y_2,z_1)\in D_x} \frac{f(x)f(y_2)f(z_1)}{(1-f(x))(1-[F(z_1)-F^{\ast}(x)])(1-[F(z_1)-F^{\ast}(y_1)])F^{\ast}(y_1)} \\
&=& \frac{1}{F^{\ast}(y_1)}\sum_{y_1<x<z_2} \frac{f(x)}{1-f(x)} \sum_{x<z_1<z_2} \frac{f(z_1)}{(1-[F(z_1)-F^{\ast}(x)])(1-[F(z_1)-F^{\ast}(y_1)])}\\
&& \sum_{y_2<y_1} f(y_2)   \\
&=& \sum_{y_1<x<z_2}\frac{f(x)}{1-f(x)} \sum_{x<z_1<z_2} \frac{f(z_1)}{(1-[F(z_1)-F^{\ast}(x)])(1-[F(z_1)-F^{\ast}(y_1)])}.
\end{eqnarray*}

\Bin From now, we use the following decomposition 

\begin{eqnarray*}
&&\frac{1}{(1-[F(z_1)-F^{\ast}(x)])(1-[F(z_1)-F^{\ast}(y_1)])} \\
&&= \frac{1}{F^{\ast}(x)-F^{\ast}(y_1)}\biggr(\frac{1}{1-[F(z_1)-F^{\ast}(y_1)]} - \frac{1}{1-[F(z_1)-F^{\ast}(x)]}\biggr)
\end{eqnarray*}

\Bin to have

\begin{eqnarray*}
B_3 &=& \sum_{y_1<x<z_2}\frac{f(x)}{(1-f(x))(F^{\ast}(x)-F^{\ast}(y_1))} \sum_{x<z_1<z_2} \frac{f(z_1)}{1-[F(z_1)-F^{\ast}(y_1)])} \\
&-& \sum_{y_1<x<z_2}\frac{f(x)}{(1-f(x))(F^{\ast}(x)-F^{\ast}(y_1))} \sum_{x<z_1<z_2} \frac{f(z_1)}{1-[F(z_1)-F^{\ast}(x)])}.
\end{eqnarray*}

\begin{eqnarray*}
B_4 &=&  \sum_{x,y_2,z_1,\ (x,y_2,z_1)\in D_x} \frac{f(x)f(y_2)f(z_1)}{(1-f(x))(1-[F(x)-F^{\ast}(y_1)])(1-[F(z_1)-F^{\ast}(y_1)])(1-F(z_1))} \\
&=& \sum_{y_1<x<z_2} \frac{f(x)}{(1-f(x))(1-[F(x)-F^{\ast}(y_1)])} \sum_{x<z_1<z_2} \frac{f(z_1)}{(1-[F(z_1)-F^{\ast}(y_1)])(1-F(z_1))} \\
&\times& \sum_{y_2<y_1} f(y_2)   \\
&=& \sum_{y_1<x<z_2} \frac{f(x)}{(1-f(x))(1-[F(x)-F^{\ast}(y_1)])} \sum_{x<z_1<z_2} \frac{f(z_1)F^{\ast}(y_1)}{(1-[F(z_1)-F^{\ast}(y_1)])(1-F(z_1))}.   
\end{eqnarray*}

\Bin From now, we use the following decomposition 

\begin{eqnarray*}
&&\frac{F^{\ast}(y_1)}{(1-[F(z_1)-F^{\ast}(y_1)])(1-F(z_1))} \\
&&=\frac{1}{1-F(z_1)} - \frac{1}{1-[F(z_1)-F^{\ast}(y_1)]}
\end{eqnarray*}

\Bin to have

\begin{eqnarray*}
B_4 &=& \sum_{y_1<x<z_2} \frac{f(x)}{(1-f(x))(1-[F(x)-F^{\ast}(y_1)])} \sum_{x<z_1<z_2} \frac{f(z_1)}{1-F(z_1)} \\
&-& \sum_{y_1<x<z_2} \frac{f(x)}{(1-f(x))(1-[F(x)-F^{\ast}(y_1)])} \sum_{x<z_1<z_2} \frac{f(z_1)}{1-[F(z_1)-F^{\ast}(y_1)]}.
\end{eqnarray*}

\begin{eqnarray*}
B_5 &=&  \sum_{x,y_2,z_1,\ (x,y_2,z_1)\in D_x} \frac{f(x)f(y_2)f(z_1)}{(1-f(x))(1-[F(x)-F^{\ast}(y_1)])(1-[F(z_1)-F^{\ast}(y_1)])F^{\ast}(y_1)} \\
&=& \frac{1}{F^{\ast}(y_1)}\sum_{y_1<x<z_2} \frac{f(x)}{(1-f(x))(1-[F(x)-F^{\ast}(y_1)])} \sum_{x<z_1<z_2} \frac{f(z_1)}{(1-[F(z_1)-F^{\ast}(y_1)])}\\
&& \sum_{y_2<y_1} f(y_2)   \\
&=& \sum_{y_1<x<z_2} \frac{f(x)}{(1-f(x))(1-[F(x)-F^{\ast}(y_1)])} \sum_{x<z_1<z_2} \frac{f(z_1)}{1-[F(z_1)-F^{\ast}(y_1)]}.   
\end{eqnarray*}

\begin{eqnarray*}
B_6 &=&  \sum_{x,y_2,z_1,\ (x,y_2,z_1)\in D_x} \frac{f(x)f(y_2)f(z_1)}{(1-f(x))(1-[F(x)-F^{\ast}(y_1)])(1-F(x))(1-F(z_1))} \\
&=& \sum_{y_1<x<z_2} \frac{f(x)}{(1-f(x))(1-[F(x)-F^{\ast}(y_1)])(1-F(x))} \sum_{x<z_1<z_2} \frac{f(z_1)}{1-F(z_1)} \sum_{y_2<y_1} f(y_2) \\   
&=& F^{\ast}(y_1)\sum_{y_1<x<z_2} \frac{f(x)}{(1-f(x))(1-[F(x)-F^{\ast}(y_1)])(1-F(x))} \sum_{x<z_1<z_2} \frac{f(z_1)}{1-F(z_1)}.
\end{eqnarray*}

\Bin Hence by regrouping all terms together, we will have

\begin{eqnarray*}
&& B_1 + B_2 + B_3 + B_4 + B_5 + B_6 \\
&& = \sum_{y_1<x<z_2} \frac{f(x)}{F^{\ast}(x)(1-F(x))} \sum_{x<z_1<z_2} \frac{f(z_1)}{1-F(z_1)}
\end{eqnarray*}

\Bin and this finishes the proof. $\blacksquare$ \\

\begin{corollary} \label{low_upp_d_03}

The \textit{pdf} of $\left(Y_{(3)},Y^{(2)}\right)$ is then given by

$$
f_{(Y_{(3)},Y^{(2)})}(y_2,z_1) = 1_{(y_2<z_1)}f(y_2)f(z_1) \sum_{y_2<x<z_1} \frac{f(x)}{F^{\ast}(x)(1-F(x))} \sum_{y_2<y_1<x} \frac{f(y_1)}{F^{\ast}(y_1)}.
$$

\end{corollary}

\Bin \textbf{Proof of Corollary \ref{low_upp_d_03}}. We have 

\begin{eqnarray*}
f_{(Y_{(3)},Y^{(2)})}(y_2,z_1) &=& \sum_{x} \sum_{y_1} \sum_{z_2} f_{Z_3}(x,y_1,y_2,z_1,z_2) \\
&=:& 1_{(y_2<z_1)} f(y_2)f(z_1)\left(C_1 + C_2 + C_3 + C_4 + C_5 + C_6\right).
\end{eqnarray*}

\Bin Let us denote by, for $(y_2,z_1)$ fixed,  $D_{(y_2,z_1)} = (y_2<y_1<x<z_1<z_2) =: D_x$ as the section of the domain of $Z_3$ at $(y_2,z_1)$ . So

\begin{eqnarray*}
C_1 &=&  \sum_{x,y_1,z_2,\ (x,y_1,z_2)\in D_x} \frac{f(x)f(y_1)f(z_2)}{(1-f(x))(1-[F(z_1)-F^{\ast}(x)])F^{\ast}(x)F^{\ast}(y_1)} \\
&=& \sum_{y_2<x<z_1} \frac{f(x)}{(1-f(x))(1-[F(z_1)-F^{\ast}(x)])F^{\ast}(x)} \sum_{y_2<y_1<x} \frac{f(y_1)}{F^{\ast}(y_1)} \sum_{z_2>z_1} f(z_2)  \\
&=& (1-F(z_1))\sum_{y_2<x<z_1} \frac{f(x)}{(1-f(x))(1-[F(z_1)-F^{\ast}(x)])F^{\ast}(x)} \sum_{y_2<y_1<x} \frac{f(y_1)}{F^{\ast}(y_1)}.
\end{eqnarray*}

\begin{eqnarray*}
C_2 &=&  \sum_{x,y_1,z_2,\ (x,y_1,z_2)\in D_x} \frac{f(x)f(y_1)f(z_2)}{(1-f(x))(1-[F(z_1)-F^{\ast}(x)])(1-[F(z_1)-F^{\ast}(y_1)])(1-F(z_1))} \\
&=& \frac{1}{1-F(z_1)}\sum_{y_2<x<z_1} \frac{f(x)}{(1-f(x))(1-[F(z_1)-F^{\ast}(x)])} \sum_{y_2<y_1<x} \frac{f(y_1)}{1-[F(z_1)-F^{\ast}(y_1)]} \sum_{z_2>z_1} f(z_2)  \\
&=& \sum_{y_2<x<z_1} \frac{f(x)}{(1-f(x))(1-[F(z_1)-F^{\ast}(x)])} \sum_{y_2<y_1<x} \frac{f(y_1)}{1-[F(z_1)-F^{\ast}(y_1)]}.
\end{eqnarray*}

\begin{eqnarray*}
C_3 &=&  \sum_{x,y_1,z_2,\ (x,y_1,z_2)\in D_x} \frac{f(x)f(y_1)f(z_2)}{(1-f(x))(1-[F(z_1)-F^{\ast}(x)])(1-[F(z_1)-F^{\ast}(y_1)])F^{\ast}(y_1)} \\
&=& \sum_{y_2<x<z_1} \frac{f(x)}{(1-f(x))(1-[F(z_1)-F^{\ast}(x)])} \sum_{y_2<y_1<x} \frac{f(y_1)}{(1-[F(z_1)-F^{\ast}(y_1)])F^{\ast}(y_1)} \sum_{z_2>z_1} f(z_2)  \\
&=& \sum_{y_2<x<z_1} \frac{f(x)}{(1-f(x))(1-[F(z_1)-F^{\ast}(x)])} \sum_{y_2<y_1<x} \frac{f(y_1)(1-F(z_1))}{(1-[F(z_1)-F^{\ast}(y_1)])F^{\ast}(y_1)}.
\end{eqnarray*}

\Bin From now, we use the following decomposition

\begin{eqnarray*}
&&\frac{1-F(z_1)}{(1-[F(z_1)-F^{\ast}(y_1)])F^{\ast}(y_1)} \\
&&= \frac{1}{F^{\ast}(y_1)} - \frac{1}{1-[F(z_1)-F^{\ast}(y_1)]}
\end{eqnarray*}

\Bin to have

\begin{eqnarray*}
C_3 &=& \sum_{y_2<x<z_1} \frac{f(x)}{(1-f(x))(1-[F(z_1)-F^{\ast}(x)])} \sum_{y_2<y_1<x} \frac{f(y_1)}{F^{\ast}(y_1)} \\
&-& \sum_{y_2<x<z_1} \frac{f(x)}{(1-f(x))(1-[F(z_1)-F^{\ast}(x)])} \sum_{y_2<y_1<x} \frac{f(y_1)}{1-[F(z_1)-F^{\ast}(y_1)]}.
\end{eqnarray*}

\begin{eqnarray*}
C_4 &=& \sum_{x,y_1,z_2,\ (x,y_1,z_2)\in D_x} \frac{f(x)f(y_1)f(z_2)}{(1-f(x))(1-[F(x)-F^{\ast}(y_1)])(1-[F(z_1)-F^{\ast}(y_1)])(1-F(z_1))} \\
&=& \frac{1}{(1-F(z_1))}\sum_{y_2<x<z_1} \frac{f(x)}{1-f(x)} \sum_{y_2<y_1<x} \frac{f(y_1)}{(1-[F(x)-F^{\ast}(y_1)])(1-[F(z_1)-F^{\ast}(y_1)])} \\
&&\sum_{z_2>z_1} f(z_2)   \\
&=& \sum_{y_2<x<z_1} \frac{f(x)}{1-f(x)} \sum_{y_2<y_1<x} \frac{f(y_1)}{(1-[F(x)-F^{\ast}(y_1)])(1-[F(z_1)-F^{\ast}(y_1)])}. 
\end{eqnarray*}

\Bin From now, we use the following decomposition

\begin{eqnarray*}
&&\frac{1}{((1-[F(x)-F^{\ast}(y_1)])(1-[F(z_1)-F^{\ast}(y_1)])} \\
&&=\frac{1}{F(z_1)-F(x)}\biggr(\frac{1}{1-[F(z_1)-F^{\ast}(y_1)]} - \frac{1}{1-[F(x)-F^{\ast}(y_1)]}\biggr)
\end{eqnarray*}

\Bin to have

\begin{eqnarray*}
C_4 &=& \sum_{y_2<x<z_1} \frac{f(x)}{(1-f(x))(F(z_1)-F(x))} \sum_{y_2<y_1<x} \frac{f(y_1)}{1-[F(z_1)-F^{\ast}(y_1)]}  \\  
&-& \sum_{y_2<x<z_1} \frac{f(x)}{(1-f(x))(F(z_1)-F(x))} \sum_{y_2<y_1<x} \frac{f(y_1)}{1-[F(x)-F^{\ast}(y_1)]}
\end{eqnarray*}

\begin{eqnarray*}
C_5 &=&  \sum_{x,y_1,z_2,\ (x,y_1,z_2)\in D_x} \frac{f(x)f(y_1)f(z_2)}{(1-f(x))(1-[F(x)-F^{\ast}(y_1)])(1-[F(z_1)-F^{\ast}(y_1)])F^{\ast}(y_1)} \\
&=& \sum_{y_2<x<z_1} \frac{f(x)}{1-f(x)} \sum_{y_2<y_1<x} \frac{f(y_1)}{(1-[F(x)-F^{\ast}(y_1)])(1-[F(z_1)-F^{\ast}(y_1)])F^{\ast}(y_1)}\\
&& \sum_{z_2>z_1} f(z_2)   \\
&=& (1-F(z_1))\sum_{y_2<x<z_1} \frac{f(x)}{1-f(x)} \\
&&\sum_{y_2<y_1<x} \frac{f(y_1)}{(1-[F(x)-F^{\ast}(y_1)])(1-[F(z_1)-F^{\ast}(y_1)])F^{\ast}(y_1)}.
\end{eqnarray*}

\Bin From now, we use the following decomposition 

\begin{eqnarray*}
&&\frac{1}{(1-[F(x)-F^{\ast}(y_1)])(1-[F(z_1)-F^{\ast}(y_1)])F^{\ast}(y_1)} \\
&&=\frac{1}{(F(z_1)-F(x))(1-F(z_1))}\biggr(\frac{1}{F^{\ast}(y_1)} - \frac{1}{1-[F(z_1)-F^{\ast}(y_1)]}\biggr) \\
&&-\frac{1}{(F(z_1)-F(x))(1-F(x))}\biggr(\frac{1}{F^{\ast}(y_1)} - \frac{1}{1-[F(x)-F^{\ast}(y_1)]}\biggr)
\end{eqnarray*}

\Bin to get

\begin{eqnarray*}
C_5 &=& \sum_{y_2<x<z_1} \frac{f(x)}{(1-f(x))(F(z_1)-F(x))} \sum_{y_2<y_1<x} \frac{f(y_1)}{F^{\ast}(y_1)} \\
&-& \sum_{y_2<x<z_1} \frac{f(x)}{(1-f(x))(F(z_1)-F(x))} \sum_{y_2<y_1<x} \frac{f(y_1)}{1-[F(z_1)-F^{\ast}(y_1)]} \\
&-& (1-F(z_1))\sum_{y_2<x<z_1} \frac{f(x)}{(1-f(x))(F(z_1)-F(x))(1-F(x))} \sum_{y_2<y_1<x} \frac{f(y_1)}{F^{\ast}(y_1)} \\
&+& (1-F(z_1))\sum_{y_2<x<z_1} \frac{f(x)}{(1-f(x))(F(z_1)-F(x))(1-F(x))} \sum_{y_2<y_1<x} \frac{f(y_1)}{1-[F(x)-F^{\ast}(y_1)]}.
\end{eqnarray*}

\begin{eqnarray*}
C_6 &=&  \sum_{x,y_1,z_2,\ (x,y_1,z_2)\in D_x} \frac{f(x)f(y_1)f(z_2)}{(1-f(x))(1-[F(x)-F^{\ast}(y_1)])(1-F(x))(1-F(z_1))} \\
&=& \frac{1}{1-F(z_1)}\sum_{y_2<x<z_1} \frac{f(x)}{(1-f(x))(1-F(x))} \sum_{y_2<y_1<x} \frac{f(y_1)}{1-[F(x)-F^{\ast}(y_1)]} \sum_{z_2>z_1} f(z_2)  \\
&=& \sum_{y_2<x<z_1} \frac{f(x)}{(1-f(x))(1-F(x))} \sum_{y_2<y_1<x} \frac{f(y_1)}{1-[F(x)-F^{\ast}(y_1)]}.
\end{eqnarray*}

\Bin Hence by regrouping all terms together, we will have

\begin{eqnarray*}
&& C_1 + C_2 + C_3 + C_4 + C_5 + C_6 \\
&& = \sum_{y_2<x<z_1} \frac{f(x)}{F^{\ast}(x)(1-F(x))} \sum_{y_2<y_1<x} \frac{f(y_1)}{F^{\ast}(y_1)}
\end{eqnarray*}

\Bin and this finishes the proof. $\blacksquare$ \\

\begin{corollary} \label{low_upp_d_04}

The \textit{pdf} of $\left(Y_{(3)},Y^{(3)}\right)$ is then given by

\begin{eqnarray*}
&&f_{(Y_{(3)},Y^{(3)})}(y_2,z_2)\\
&&= 1_{(y_2<z_2)}f(y_2)f(z_2) \sum_{y_2<x<z_2} \frac{f(x)}{F^{\ast}(x)(1-F(x))} \sum_{y_2<y_1<x} \frac{f(y_1)}{F^{\ast}(y_1)}\sum_{x<z_1<z_2} \frac{f(z_1)}{1-F(z_1)}.
\end{eqnarray*}
\end{corollary}

\Bin \textbf{Proof of Corollary \ref{low_upp_d_04}}. We have 

\begin{eqnarray*}
f_{(Y_{(3)},Y^{(3)})}(y_2,z_2) &=& \sum_{x} \sum_{y_1} \sum_{z_1} f_{Z_3}(x,y_1,y_2,z_1,z_2) \\
&=:& 1_{(y_2<z_2)} f(y_2)f(z_2)\left(D_1 + D_2 + D_3 + D_4 + D_5 + D_6\right).
\end{eqnarray*}

\Bin Let us denote by, for $(y_2,z_2)$ fixed,  $D_{(y_2,z_2)} = (y_2<y_1<x<z_1<z_2) =: D_x$ as the section of the domain of $Z_3$ at $(y_2,z_2)$ . So

\begin{eqnarray*}
D_1 &=&  \sum_{x,y_1,z_1,\ (x,y_1,z_1)\in D_x} \frac{f(x)f(y_1)f(z_1)}{(1-f(x))(1-[F(z_1)-F^{\ast}(x)])F^{\ast}(x)F^{\ast}(y_1)} \\
&=& \sum_{y_2<x<z_2} \frac{f(x)}{(1-f(x))F^{\ast}(x)} \sum_{x<z_1<z_2} \frac{f(z_1)}{1-[F(z_1)-F^{\ast}(x)]}\sum_{y_2<y_1<x}\frac{f(y_1)}{F^{\ast}(y_1)}. 
\end{eqnarray*}

\begin{eqnarray*}
&&D_2\\
&&=  \sum_{x,y_1,z_1,\ (x,y_1,z_1)\in D_x} \frac{f(x)f(y_1)f(z_1)}{(1-f(x))(1-[F(z_1)-F^{\ast}(x)])(1-[F(z_1)-F^{\ast}(y_1)])(1-F(z_1))} \\
&&= \sum_{y_2<x<z_2} \frac{f(x)}{1-f(x)} \sum_{x<z_1<z_2} \frac{f(z_1)}{(1-[F(z_1)-F^{\ast}(x)])(1-F(z_1))}\sum_{y_2<y_1<x}\frac{f(y_1)}{1-[F(z_1)-F^{\ast}(y_1)]}. 
\end{eqnarray*}

\begin{eqnarray*}
&&D_3\\
&&=  \sum_{x,y_1,z_1,\ (x,y_1,z_1)\in D_x} \frac{f(x)f(y_1)f(z_1)}{(1-f(x))(1-[F(z_1)-F^{\ast}(x)])(1-[F(z_1)-F^{\ast}(y_1)])F^{\ast}(y_1)} \\
&&= \sum_{y_2<x<z_2} \frac{f(x)}{1-f(x)} \sum_{x<z_1<z_2} \frac{f(z_1)}{1-[F(z_1)-F^{\ast}(x)]}\sum_{y_2<y_1<x} \frac{f(y_1)}{(1-[F(z_1)-F^{\ast}(y_1)])F^{\ast}(y_1)}. 
\end{eqnarray*}

\Bin From now, we use the following decomposition 

\begin{eqnarray*}
&&\frac{1}{(1-[F(z_1)-F^{\ast}(y_1)])F^{\ast}(y_1)} \\
&&=\frac{1}{1-F(z_1)}\biggr(\frac{1}{F^{\ast}(y_1)}-\frac{1}{1-[F(z_1)-F^{\ast}(y_1)]}\biggr)
\end{eqnarray*}

\Bin to have
 
\begin{eqnarray*}
D_3 &=& \sum_{y_2<x<z_2} \frac{f(x)}{1-f(x)} \sum_{x<z_1<z_2} \frac{f(z_1)}{(1-F(z_1))(1-[F(z_1)-F^{\ast}(x]}\sum_{y_2<y_1<x}\frac{f(y_1)}{F^{\ast}(y_1)} \\
&-& \sum_{y_2<x<z_2} \frac{f(x)}{1-f(x)} \sum_{x<z_1<z_2} \frac{f(z_1)}{(1-F(z_1))(1-[F(z_1)-F^{\ast}(x)])}\sum_{y_2<y_1<x}\frac{f(y_1)}{1-[F(z_1)-F^{\ast}(y_1)]}. 
\end{eqnarray*}

\begin{eqnarray*}
D_4 &=&  \sum_{x,y_1,z_1,\ (x,y_1,z_1)\in D_x} \frac{f(x)f(y_1)f(z_1)}{(1-f(x))(1-[F(x)-F^{\ast}(y_1)])(1-[F(z_1)-F^{\ast}(y_1)])(1-F(z_1))} \\
&=& \sum_{y_2<x<z_2} \frac{f(x)}{1-f(x)} \sum_{x<z_1<z_2} \frac{f(z_1)}{1-F(z_1)}\sum_{y_2<y_1<x} \frac{f(y_1)}{(1-[F(x)-F^{\ast}(y_1)])(1-[F(z_1)-F^{\ast}(y_1)])}. 
\end{eqnarray*}

\Bin From now, we use the following decomposition

\begin{eqnarray*}
&&\frac{1}{(1-[F(x)-F^{\ast}(y_1)])(1-[F(z_1)-F^{\ast}(y_1)])} \\
&&=\frac{1}{F(z_1)-F(x)}\biggr(\frac{1}{1-[F(z_1)-F^{\ast}(y_1)]} - \frac{1}{1-[F(x)-F^{\ast}(y_1)]}\biggr) 
\end{eqnarray*}
 
\Bin to have 

\begin{eqnarray*}
D_4 &=& \sum_{y_2<x<z_2} \frac{f(x)}{1-f(x)} \sum_{x<z_1<z_2} \frac{f(z_1)}{(F(z_1)-F(x))(1-F(z_1))}\sum_{y_2<y_1<x} \frac{f(y_1)}{1-[F(z_1)-F^{\ast}(y_1)]} \\
&-& \sum_{y_2<x<z_2} \frac{f(x)}{1-f(x)} \sum_{x<z_1<z_2} \frac{f(z_1)}{(F(z_1)-F(x))(1-F(z_1))}\sum_{y_2<y_1<x} \frac{f(y_1)}{1-[F(x)-F^{\ast}(y_1)]}.
\end{eqnarray*}

\begin{eqnarray*}
D_5 &=&  \sum_{x,y_1,z_1,\ (x,y_1,z_1)\in D_x} \frac{f(x)f(y_1)f(z_1)}{(1-f(x))(1-[F(x)-F^{\ast}(y_1)])(1-[F(z_1)-F^{\ast}(y_1)])F^{\ast}(y_1)} \\
&=& \sum_{y_2<x<z_2} \frac{f(x)}{1-f(x)} \sum_{x<z_1<z_2} f(z_1)\sum_{y_2<y_1<x} \frac{f(y_1)}{(1-[F(x)-F^{\ast}(y_1)])(1-[F(z_1)-F^{\ast}(y_1)])F^{\ast}(y_1)}. 
\end{eqnarray*}

\Bin From now, we use the following decomposition 

\begin{eqnarray*}
&&\frac{1}{(1-[F(x)-F^{\ast}(y-1)])(1-[F(z_1)-F^{\ast}(y_1)])F^{\ast}(y_1)} \\
&&=\frac{1}{(F(z_1)-F(x))(1-F(z_1))}\biggr(\frac{1}{F^{\ast}(y_1)} - \frac{1}{1-[F(z_1)-F^{\ast}(y_1)]}\biggr) \\
&&-\frac{1}{(F(z_1)-F(x))(1-F(x))}\biggr(\frac{1}{F^{\ast}(y_1)} - \frac{1}{1-[F(x)-F^{\ast}(y_1)]}\biggr)
\end{eqnarray*}

\Bin to get

\begin{eqnarray*}
&&D_5\\
&&= \sum_{y_2<x<z_2} \frac{f(x)}{1-f(x)} \sum_{x<z_1<z_2} \frac{f(z_1)}{(F(z_1)-F(x))(1-F(z_1))}\sum_{y_2<y_1<x} \frac{f(y_1)}{F^{\ast}(y_1)} \\
&&- \sum_{y_2<x<z_2} \frac{f(x)}{1-f(x)} \sum_{x<z_1<z_2} \frac{f(z_1)}{(F(z_1)-F(x))(1-F(z_1))}\sum_{y_2<y_1<x} \frac{f(y_1)}{1-[F(z_1)-F^{\ast}(y_1)]} \\
&&- \sum_{y_2<x<z_2} \frac{f(x)}{(1-f(x))(1-F(x))} \sum_{x<z_1<z_2} \frac{f(z_1)}{F(z_1)-F(x)}\sum_{y_2<y_1<x} \frac{f(y_1)}{F^{\ast}(y_1)} \\
&&+ \sum_{y_2<x<z_2} \frac{f(x)}{(1-f(x))(1-F(x))} \sum_{x<z_1<z_2} \frac{f(z_1)}{F(z_1)-F(x)}\sum_{y_2<y_1<x} \frac{f(y_1)}{1-[F(x)-F^{\ast}(y_1)]}.
\end{eqnarray*}

\begin{eqnarray*}
D_6 &=&  \sum_{x,y_1,z_1,\ (x,y_1,z_1)\in D_x} \frac{f(x)f(y_1)f(z_1)}{(1-f(x))(1-[F(x)-F^{\ast}(y_1)])(1-F(x))(1-F(z_1))} \\
&=& \sum_{y_2<x<z_2} \frac{f(x)}{(1-f(x))(1-F(x))} \sum_{x<z_1<z_2} \frac{f(z_1)}{1-F(z_1)}\sum_{y_2<y_1<x} \frac{f(y_1)}{1-[F(x)-F^{\ast}(y_1)]}. 
\end{eqnarray*}

\Bin Hence by doing simple computations, we will have

$$
D_1 + D_2 + D_3 = \sum_{y_2<x<z_2} \frac{f(x)}{(1-f(x))F^{\ast}(x)} \sum_{x<z_1<z_2} \frac{f(z_1)}{1-F(z_1)}\sum_{y_2<y_1<x} \frac{f(y_1)}{F^{\ast}(y_1)}
$$

\Bin and

$$
D_4 + D_5 + D_6 = \sum_{y_2<x<z_2} \frac{f(x)}{(1-f(x))(1-F(x))} \sum_{x<z_1<z_2} \frac{f(z_1)}{1-F(z_1)}\sum_{y_2<y_1<x} \frac{f(y_1)}{F^{\ast}(y_1)}
$$

\Bin and the proof is over by putting all the terms together. $\blacksquare$\\

\section*{Conclusions and perspectives}

\Ni Up to three records, we have been able to establish the joint simultaneously lower and upper strong records in both frames of absolutely continuous and discrete
random variables. We used the \textit{pdf}'s to rediscover the joint \textit{pdf} for upper records only or for lower records only. The main lesson is that: for $n \in \{2, \ 3\}$, the density of the \textit{sjlu} record values \textit{pdf} or \textit{mpf} (both denoted $f(\circ)$ below) in $(y_1, \cdots,y_n, x_1,\cdots,x_n)$ is product of \\

\Ni (a) a main part of the form 

$$
\biggr(\prod_{J=1}^{N} f(x_j)f(y_j)\biggr)\ 1_{(y_1<\cdots<y_n<x_1<\cdots<x_n)}
$$

\Ni (b) by a function which is the addition of 

$$
\frac{(2(n-1))!}{((n-1)!)^2}
$$

\Ni functions $h_j$, each of them being the product of $(2n-1)$ factors and being a function of terms $f(y_j)$, $f(x_j)$, $F(x_j)$, $F(y_j)$, $1-F(x_j)$, $1-F(y_j)$,
$F(x_h)-F(y_k)$. Each function correspond to a permutation with repetitions of $(2(n-1))$ objects with two subgroups of indistinguishable objects of size $(n-1)$ each. Such a finding will be exploited in the study devoted to the general law.\\

\Bin For $n \in \{2, \ 3\}$, we extracted all the marginal laws from the \textit{sjlu} record values probability laws and then, rediscovered  the laws of pairs of upper record values and pairs of lower record values, accordingly to known results in the field.\\

\end{document}